\documentclass[11pt,reqno]{amsart}
\usepackage[tmargin=1.5in,bmargin=1.4in,rmargin=1.4in,lmargin=1.4in]{geometry}

\usepackage{amsmath,amsthm,amssymb,mathrsfs}
\usepackage[breaklinks=true]{hyperref}

\theoremstyle{plain}
\newtheorem{theorem}{Theorem}[section]
\newtheorem{proposition}[theorem]{Proposition}

\theoremstyle{definition}
\newtheorem{definition}[theorem]{Definition}
\newtheorem{remark}[theorem]{Remark}

\numberwithin{equation}{section}

\newcommand{\ba}{\mathbf{a}}
\newcommand{\balpha}{{\boldsymbol{\alpha}}}
\newcommand{\bb}{\mathbf{b}}
\newcommand{\bone}{{\boldsymbol{1}}}
\newcommand{\bX}{\mathbf{X}}
\newcommand{\C}{\mathbb{C}}
\newcommand{\cB}{\mathcal{B}}
\newcommand{\cN}{\mathcal{N}}
\newcommand{\diag}{\mathop{\mathrm{diag}}}
\newcommand{\E}{\mathbb{E}}
\newcommand{\I}{\mathrm{i}}
\newcommand{\intd}{\,\rd}
\newcommand{\PF}{\mathrm{PF}}
\newcommand{\R}{\mathbb{R}}
\renewcommand{\Re}{\mathop{\mathrm{Re}}}
\newcommand{\rd}{\mathrm{d}}
\newcommand{\sFTN}{\mathscr{F}^{\TN}}
\newcommand{\sFTP}{\mathscr{F}^{\TP}}
\newcommand{\TN}{\mathrm{TN}}
\newcommand{\TP}{\mathrm{TP}}
\newcommand{\tr}{\mathop{\mathrm{tr}}\nolimits}
\newcommand{\Z}{\mathbb{Z}}

\makeatletter
\def\blfootnote{\gdef\@thefnmark{}\@footnotetext}
\makeatother

\begin{document}

\title[Preservers of TP kernels and PF functions]%
{Preservers of totally positive kernels\\and P\'olya frequency functions}

\author{Alexander Belton}
\address[A.~Belton]{Department of Mathematics and Statistics, Lancaster University, Lancaster, UK}
\email{\tt a.belton@lancaster.ac.uk}

\author{Dominique Guillot}
\address[D.~Guillot]{University of Delaware, Newark, DE, USA}
\email{\tt dguillot@udel.edu}

\author{Apoorva Khare}
\address[A.~Khare]{Department of Mathematics, Indian Institute of
Science, Bangalore, India; and Analysis and Probability Research Group,
Bangalore, India}
\email{\tt khare@iisc.ac.in}

\author{Mihai Putinar}
\address[M.~Putinar]{University of California at Santa Barbara, CA, USA;
and Newcastle University, Newcastle upon Tyne, UK} 
\email{\tt mputinar@math.ucsb.edu, mihai.putinar@ncl.ac.uk}

\blfootnote{Received March 11, 2022; revised May 18, 2022.}
\date{}



\begin{abstract}
Fractional powers and polynomial maps preserving structured
totally positive matrices, one-sided P\'olya frequency functions, or
totally positive kernels are treated from a unifying
perspective. Besides the stark rigidity of the polynomial transforms,
we unveil an ubiquitous separation between discrete and continuous
spectra of such inner fractional powers. Classical works of
Schoenberg, Karlin, Hirschman, and Widder are completed by our
classification. Concepts of probability theory, multivariate statistics,
and group representation theory naturally enter into the picture.
\end{abstract}

\maketitle

\section{Overview}

The purpose of this note is to announce some of the results in a
sequence of three closely related recent
papers~\cite{BGKP-TN,Khare2020,BGKP-HW}
which study P\'olya frequency functions and the post-composition
transforms that preserve them and other classes of totally positive
kernels. Specifically, we complete results in total
positivity due to Karlin's \textit{Transactions}~\cite{KarlinTAMS}
and Schoenberg's \textit{Annals}~\cite{Schoenberg55}
by providing extensions and converses. These concern
fractional powers of totally non-negative Toeplitz kernels and
of their Laplace transforms. We isolate a spectral-threshold
phenomenon in total positivity which is similar to the structure of
Berezin--Gindikin--Wallach-type sets, but which was discovered in
particular situations prior to the work of these mathematicians.

Next, we focus on powers and other polynomial preservers of more
general non-smooth one-sided P\'olya frequency (PF) functions. The
underlying theme is probabilistic: one-sided PF functions are the
density functions of linear combinations $\sum_j \alpha_j X_j$ of
independent standard exponential random variables: for example,
Karlin's 1964 kernel is the density of $X_1 + X_2$.
The densities $\Lambda_\balpha$ of finite linear combinations
$\sum_{j=1}^m \alpha_j X_j$ were studied by Hirschman and Widder
\cite{HW49}, and we prove that, if the
coefficients $\balpha = ( \alpha_j ) \in ( 0, \infty )^m$ lie
outside a null set then the only polynomials~$p$ with the property that
$p \circ \Lambda_\balpha$ is a PF function, are homotheties
$p(x) = c x$, where $c > 0$.  A previously unexplored
connection between recovering
$\balpha$ from the moments of $\Lambda_\balpha$ and the Jacobi--Trudi
identity from symmetric function theory comes as a bonus.

Finally, we study general one-sided PF functions, which are the densities
of countable sums of exponential random variables, as well as other
classes of totally non-negative kernels, including general PF functions
and sequences. For these classes of maps, we show there are very few
preservers beyond positive homotheties $p( x ) = c x$. The culmination of
this line of inquiry is the characterization of preservers of
total positivity and total non-negativity for kernels on $X \times Y$,
where $X$ and $Y$ are arbitrary totally ordered domains.

In addition we provide a small correction to Schoenberg's classification
of discontinuous PF functions~\cite{Schoenberg51}.

\section{The pervasive nature of P\'olya frequency functions}

\subsection{Total positivity and P\'olya frequency functions}

Given totally ordered sets $X$ and $Y$, a kernel
$K : X \times Y \to \R$ is \emph{totally positive} if, for all integers
$n \geq 1$ and all choices of $x_1 < \cdots < x_n$ in $X$ and $y_1 <
\cdots < y_n$ in $Y$, the determinant $\det ( K( x_i, y_j ) )_{i, j =
1}^n$ is positive;
the kernel $K$ is \emph{totally non-negative} if these determinants are
non-negative.
We refer to such kernels, which are matrices if the domains $X$ and $Y$
are finite, as $\TP$ and $\TN$ kernels, respectively. If non-negativity
or positivity are only required to hold whenever $n \leq p$ then we speak
of $\TN_p$ or $\TP_p$ kernels, respectively.
(The monographs \cite{Karlin, Pinkus} refer to \textit{strict total
positivity} and \textit{total positivity} instead of total positivity and
total non-negativity.)

Total positivity is a long-studied and evergreen area of mathematics. For
almost a century, totally positive and totally non-negative kernels
surfaced in the most unexpected circumstances, and this trend continues
in full force today. Although this chapter of matrix analysis remains
somewhat recondite, it has reached maturity due to the dedicated efforts
of several generations of mathematicians. The foundational work
by Gantmacher and Krein \cite{GK},
the survey \cite{Ando}, the early monograph \cite{GK1}, and the more
recent publications \cite{Karlin,TP,Pinkus,FJ} offer ample references to
the fascinating history of total positivity, as well as accounts of its
many surprising applications. These include
analysis~\cite{Schoenberg50,Schoenberg51,Schoenberg55},
differential equations~\cite{Loewner55},
probability and statistics~\cite{Johnstone,Efron,Karlin,KP}, 
and interpolation theory~\cite{Curry2,SchoenbergWhitney53},
to provide a few areas and early references.
Total positivity continues to make very recent impacts in areas such as
representation theory and cluster algebras
\cite{BFZ,BZ-2,FZ-1,FZ-2,L-1,L-2,Ri}, Gabor analysis \cite{GRS,GS},
combinatorics \cite{Br1,Br2}, as well as integrable systems and positive
Grassmannians (the geometric avatar of total positivity)
\cite{KW-1,KW-2,Postnikov}.

The origins of total positivity lie in the property of diminishing
variation, which can be traced back to Descartes
(1600s, via his rule of signs~\cite{Descartes}), but concretely at
least as far back as Laguerre~\cite{Laguerre2} (1883) and
Fekete~\cite{Fe} (1912). P\'olya then coined the phrase
`variationsvermindernd', and Schoenberg showed in \cite{Schoenberg30}
(1930) that $\TP$ and $\TN$ matrices have this property of variation
diminution. There has been continuing activity for $\TN$
matrices~\cite{GK1,Johnstone,Projesh}, and also for $\TN$ kernels on
bi-infinite domains, which are the main focus of this note.

\begin{definition}
A function $\Lambda : \R \to \R$ is
\textit{totally non-negative} if
the associated Toeplitz kernel $T_\Lambda$ is
totally non-negative, where
\begin{equation}\label{ETlambda}
T_\Lambda : \R \times \R \to \R; \ ( x, y ) \mapsto \Lambda( x - y ).
\end{equation}
If, further, the function $\Lambda$ is Lebesgue integrable and non-zero
at two or more points, then $\Lambda$ is a
\textit{P\'olya frequency function}.
A function $\Lambda : \Z \to \R$ whose associated Toeplitz
kernel is totally non-negative is a \textit{P\'olya frequency sequence}.
\end{definition}

Following early work by P\'olya and Hamburger, Schoenberg has initiated
and developed the theory of P\'olya frequency functions in his landmark
paper~\cite{Schoenberg51} dated 1951 (following announcements in 1947 and
$1948$ in \textit{Proc.\ Natl.\ Acad.\ Sci.\ USA}). In the companion
work~\cite{Schoenberg50} (published a year earlier),
Schoenberg proved that, when viewed as convolution integral operators,
P\'olya frequency functions can be characterized in terms of the
variation-diminishing property. This study led to an explosion of work
in numerical analysis and approximation theory, via \emph{splines}.

\subsection{The Laguerre--P\'olya class of entire functions}

A second notable appearance of P\'olya frequency functions and sequences is within
the theory of functions of a complex variable. We henceforth refer to these classes of maps as
\emph{$\PF$ functions} and \emph{$\PF$ sequences}, respectively.

The Fourier--Laplace transform of a $\PF$ function is the reciprocal,
on a suitable domain of definition, of a \emph{Laguerre--P\'olya entire function}~\cite{ASW,Edrei,Schoenberg51}.
The natural question of characterizing locally uniform limits of sequences of polynomials with only real roots
was answered by Laguerre \cite{Laguerre1}, and completed by P\'olya \cite{Polya0}.
A related wider program of locating the zeros of entire functions was initiated about a century ago by
P\'olya and Schur \cite{Polya-Schur}. By now, the topic of Laguerre--P\'olya entire functions appears in textbooks \cite{Levin},
reflecting a century of accumulated knowledge passing through early contributions such as \cite{Grommer, Tschebotareff} and
continuing up to the present \cite{BCC}. Since its inception, Riemann's hypothesis was and remains a background theme in studies
concerning entire functions of the Laguerre--P\'olya class, and so indirectly to P\'olya frequency functions and sequences.
We cite here only two recent contributions: \cite{Grochenig,Katkova2007}.

We owe to Schoenberg the aforementioned characterization of $\PF$ functions.

\begin{theorem}[Schoenberg, \cite{Schoenberg51}]\label{Tlaplace}
Given a P\'olya frequency function $\Lambda$, its bilateral Laplace
transform
\[
\cB\{ \Lambda \}( s ) := \int_\R e^{-x s} \Lambda( x ) \intd x
\]
converges for complex $s$ in an open strip containing the imaginary axis,
and equals $1/\Psi$ on this strip, for an entire function $\Psi$ in the
Laguerre--P\'olya class. Conversely, any function $\Psi$ of the above
form agrees with the reciprocal of the bilateral Laplace transform of
some P\'olya frequency function on its strip of convergence.
\end{theorem}

In \cite{Schoenberg51}, Schoenberg then showed that a $\PF$ function
either vanishes precisely on a semi-axis, which may be open or closed,
or is non-vanishing on $\R$. These functions are termed one-sided and
two-sided $\PF$ functions, respectively. Schoenberg also showed in
\cite{Schoenberg51} that, up to linear transformations, the
reciprocals of their Laplace transforms are respectively in the first
Laguerre--P\'olya class of entire functions~\cite{Laguerre1,Polya0}
that are non-vanishing at~$0$,
\begin{equation}\label{Elp1}
\Psi(s) = C e^{\delta s} \prod_{j=1}^\infty (1 + \alpha_j s), \quad
\text{with}\ C > 0, \ \delta, \alpha_j \geq 0, \ \sum_j
\alpha_j < \infty,
\end{equation}
and the second Laguerre--P\'olya class of functions not vanishing at
the origin,
\begin{equation}\label{Elp2}
\Psi(s) = C e^{-\gamma s^2 + \delta s}
\prod_{j=1}^\infty (1 + \alpha_j s) e^{- \alpha_j s},
\quad \text{with} \ C > 0, \ \gamma \geq 0, \ \delta, \alpha_j \in \R, 
\ \sum_j \alpha_j^2 < \infty.
\end{equation}

This remarkable dictionary, established by Schoenberg, provides the building blocks of P\'olya frequency functions.
This class of functions is closed under convolution, and has Gaussian functions,
one-sided exponentials and simple fractions with poles on the real line as generators.
The exploitation of this basic observation is paramount for interpolation theory, and the
expository lectures by Schoenberg \cite{Schoenberg1973} are as fresh and informative
today as they were fifty years ago.

\subsection{Group representations}

The string of discoveries of the same class of objects does not stop at entire function theory.
The classification of characters of irreducible unitary representations of the infinite symmetric group
$S(\infty) = \cup_n S(n)$ and the infinite unitary group $U(\infty) = \cup_n U(n)$
led Thoma \cite{Thoma} and Voiculescu \cite{Voiculescu} independently to
the class of Fourier transforms of P\'olya frequency functions.
These remarkable findings in the 1960s and 1970s established the foundations of
the representation theory of ``big groups'': see \cite{OV,VK} for details.

Not unrelatedly, the computations of spherical functions pointed to orbital-integral formulae for such character functions,
as did the more general question of spectral synthesis on homogeneous spaces.
The pioneering work of Gelfand and Naimark \cite{GN} opened a whole chapter of explicit expressions linking orbital integrals to invariants of finite groups,
in the spirit of Weyl's character formula.

A typical orbital integral has the form
\begin{equation}\label{Eorbital}
f_A : H(n) \to \C; \ B \mapsto \int_\Omega e^{\I \tr(B M)} \mu( \rd M ),
\end{equation}
where $A \in H(n)$ is a positive semi-definite $n \times n$ complex matrix, $\Omega$ denotes its orbit under conjugation by the unitary group $U(n)$,
and $\mu$ is a $U(n)$-invariant measure carried by $\Omega$.
In view of the invariance of the trace under cyclic permutations, $f_A$ is invariant under unitary conjugation:
\[
f_A( U B U^\ast ) = f_A( B ) \qquad \text{for all } U \in U(n) \text{ and } B \in H(n).
\]
In particular $f_A( B )$ depends only on the eigenvalues of $B$ and is a symmetric function of these eigenvalues.
Furthermore, being a Fourier transform, the function $f_A$ is positive definite.

To bring P\'olya frequency functions into view, one considers the inductive limit of such measures and functions
defined on the union $H(\infty) = \bigcup_n H(n)$ of Hermitian matrices of arbitrary size.
These spherical functions, normalized by the condition that $f( 0 ) = 1$, form a
convex set and the extremal points of this set are multiplicative, in the sense that
\[
f\bigl( \diag(b_1, b_2, \ldots, b_m) \bigr) = F( b_1 ) F( b_2 ) \cdots F( b_m )
\]
for some function of a real variable $F$. This occurs
precisely when the corresponding unitarily invariant measure $\mu$ on the union $H(\infty)$ is ergodic.
The main classification theorem of \cite[Theorem~2.9]{OV} asserts the existence of a {\it  bijective
correspondence between ergodic, unitarily invariant probability measures on
$H(\infty)$ and P\'olya frequency functions}. To be more precise, $F$
is the Fourier transform of a P\'olya frequency function associated with the
ergodic measure $\mu$. Moreover, specific invariant measures provide the building blocks
of the class of P\'olya frequency functions \cite[Corollaries~2.5 to~2.7]{OV}.

Let us consider the case when $A = \diag( a_1, \ldots, a_m )$, where $a_1$, \ldots, $a_m$ are positive,
and let $B = E_{11} = \diag( 1, 0, 0, \ldots, 0 )$. Passing lightly over the technicalities required to extend~$f_A$, and
using the symmetry $f_A( \I x B ) = f_B( \I x A )$, which follows from the tracial property, we have that
\begin{equation}\label{Eorbit}
f_A( \I x E_{11} ) = \int_{\Omega'} \exp\Bigl(-x  \sum_{j=1}^m a_j |z_j|^2 \Bigr) \sigma( \rd z ) %
\qquad ( x > 0 ),
\end{equation}
where $\Omega' = S^{2 m - 1} \cong U( m )/ U( m - 1 )$ is the unit sphere in $\C^m$
and $\sigma$ is the normalized rotationally invariant measure on the sphere.

Up to proper normalization, the spherical average (\ref{Eorbital}) is
equal to the Hirschman--Widder distribution $\Lambda_\balpha$ at the point~$x$,
where $\balpha = ( a_1^{-1}, \ldots, a_m^{-1} )$.
The right-hand side of the identity (\ref{Eorbit}) is known as a Harish-Chandra--Itzykson--Zuber integral \cite{hc,iz}.

We can go further and establish an explicit link between the
orbital integral $f_A(B)$ above and Hirschman--Widder densities.
Further details, including complete proofs, are contained in \cite[Section 5]{OV}.
Let $\ba = ( a_1, \ldots, a_m ) \in \R^m$ and $\bb = ( b_1, \ldots, b_m ) \in \C^m$
have corresponding diagonal matrices $A = \diag \ba \in H(m)$ and $B = \diag \bb$, respectively.
As in (\ref{Eorbital}), we let $f_A$ denote the characteristic function of the invariant probability
measure~$\mu$ with support~$\Omega$, where $\Omega$ is the $U(m)$-orbit of $A$ under conjugation:
\[
f_A( B ) := \int_\Omega e^{\I \tr( B M )} \mu( \rd M ) =  \int_{U(m)} e^{\I \tr( B U A U^* )} \intd U.
\]
Since $f_A$ is entire and symmetric as a function of the coordinates of $\bb$,
it admits a Taylor-series expansion that is convergent everywhere, so also a convergent expansion in terms of Schur
polynomials:
\[
f_A( \diag \bb ) = \sum_\nu c_\nu s_\nu( \bb ),
\]
where the sum runs over Young diagrams with at most $m$ rows. A
computation by Olshanski and Vershik, using characters of $U(m)$ and
change-of-bases formulas between symmetric power-sum polynomials and Schur
polynomials, provides a closed-form expression for the coefficient
$c_\nu$: see \cite[Theorem~5.1]{OV}. This strategy
appeared a few decades earlier, in explicit computations of Gel'fand and Naimark \cite{GN}, and
quite remarkably (and independently) in multivariate statistics: see James \cite{James} and
the comments in \cite{Farrell}. From here, one derives the following
expansions: see \cite[Corollaries~5.2 and~5.4]{OV}.

\begin{proposition}\label{POV}
If the tuples $\ba = ( a_1, \ldots, a_m ) \in \R^m$ and $\bb = ( b_1, \ldots, b_m ) \in \C^m$
each have distinct coordinates and $A = \diag \ba$
then the orbital integral $f_A$ is given by the Harish-Chandra--Itzykson--Zuber formula:
\[
f_A( -\I \diag \bb ) = \frac{\prod_{j = 0}^{m - 1} j!}{V( \ba ) V( \bb )}
\det\begin{pmatrix}
e^{b_1 a_1} & e^{b_1 a_2} & \cdots & e^{b_1 a_m}\\
e^{b_2 a_1} & e^{b_2 a_2} & \cdots & e^{b_2 a_m}\\
 \vdots & \vdots & \ddots & \vdots \\
e^{b_m a_1} & e^{b_m a_2} & \cdots & e^{b_m a_m}
\end{pmatrix}.
\]
If instead $B = \diag( 1, 0, \ldots, 0 ) = E_{11}$ then
\[
f_A( -\I x B ) = ( m - 1 )! \sum_{j = 0}^\infty \frac{h_j( a_1, \ldots, a_m )}{( j + m - 1 )!} x^j,
\]
where $h_j$ is the $j$th complete homogeneous symmetric polynomial.
\end{proposition}

In particular, if $a_1$, \ldots, $a_m$ are positive and distinct, and $x > 0$ then,
by the second part of Proposition~\ref{POV} and the identity (\ref{Emaclaurinhw}) below,
\[
f_A( \I x E_{11} ) = %
( m - 1 )! ( -x )^{1 - m} \sum_{n = m - 1}^\infty \frac{h_{n - m + 1}( a_1, \ldots,  a_m )}{n!} ( -x )^n
= \frac{( m - 1 )! x^{1 - m}}{a_1 \cdots a_m} \Lambda_\balpha( x ),
\]
where $\balpha = ( a_1^{-1}, \ldots,  a_m^{-1} )$.

In conclusion, the Hirschman--Widder density possesses the following integral and determinantal representations: if $x > 0$ then
\begin{align*}
\Lambda_\balpha( x ) & = \frac{a_1 \cdots a_m}{( m - 1 )!} x^{m - 1} %
\int_{S^{2 m - 1}} \exp\Bigl(-x \sum_{j=1}^m a_j |z_j|^2 \Bigr) \sigma( \rd z) \\
 & = \frac{a_1 \cdots a_m}{V( \ba )}
 \det\begin{pmatrix}
e^{-a_1 x} & e^{-a_2 x} & e^{-a_3 x} & \cdots & e^{-a_m x} \\
1 & 1 & 1 & \cdots & 1 \\
a_1 & a_2 & a_3 & \cdots & a_m\\
\vdots & \vdots & \vdots & \ddots & \vdots\\
a_1^{m - 2} & a_2^{m - 2} & a_3^{m - 2} & \cdots & a_m^{m - 2}
\end{pmatrix}.
\end{align*}
The second representation can be obtained from the first identity in
Proposition~\ref{POV} by setting $\bb = ( -x, 0, y_1, \ldots, y_{m - 2} )$ and
taking successively the $j$th partial derivative at zero with respect to
$y_j$ for $j = 1$ to $m - 2$.

\subsection{Probability theory}

The aspects of P\'olya frequency functions described above lead to some natural interpretations of a probabilistic flavor.
Here we discuss three such perspectives.

First we return to the Hadamard--Weierstrass factorizations of Laplace
transforms of P\'olya frequency functions. Suppose we are given some
continuous random variables, whose associated probability density
functions admit bilateral Laplace transforms convergent in some common
open strip containing the imaginary axis. Linear combinations of these
random variables correspond to convolutions of their densities, which in
turn correspond to products at the level of Laplace transforms. From the
perspective of (\ref{Elp1}) and (\ref{Elp2}), Schoenberg's
Theorem~\ref{Tlaplace} can be recast using countably many exponential
random variables and a Gaussian random variable.

\begin{theorem}\label{TPFprob}
Given a P\'olya frequency function $\Lambda : \R \to \R$, exactly one of
the following holds:
\begin{enumerate}
\item $\Lambda$ is discontinuous at exactly one point $x_0$, and
vanishes on one side of $x_0$. Moreover, ignoring the value at $x_0$
and up to a shift of argument and positive rescaling,
$\Lambda$ equals the density of $\alpha X$, where $\alpha \in \R
\setminus \{ 0 \}$ and $X$ is a standard exponential random variable with
mean $1$.

\item $\Lambda$ is continuous, vanishes on a semi-axis
$\pm(-\infty,x_0)$, and is non-vanishing on the interior of the complement.
In this case, up to a shift in argument and positive rescaling,
$\Lambda$ equals the density of a linear combination
$\pm\sum_j \alpha_j X_j$, where the $\alpha_j$ are positive and summable, and the $X_j$ are
independent standard exponential variables.

\item $\Lambda$ is non-vanishing on $\R$. In this case,
up to a shift in argument and positive rescaling, $\Lambda$ equals
the density of a linear combination
$\alpha Y \pm \sum_j \alpha_j X_j$, where
(a)~the series $\sum_j \alpha_j^2$ converges,
(b)~if $\alpha=0$ then there exist at least one positive and one negative
$\alpha_j$,
(c)~the $X_j$ are independent standard exponential random variables, and
(d)~$Y$ is a standard Gaussian random variable independent of the $X_j$.
\end{enumerate}
\end{theorem}

As straightforward as it may be, we were unable to find this observation in the literature. We
recorded it in a series of results and remarks in~\cite[Section~2.4]{BGKP-HW}. In particular,
the observation about the indeterminateness of a PF function at its point
of discontinuity revealed a gap in the literature,
which is now clarified in~\cite{Khare2020}; see Theorem~\ref{Tdiscont} and the
discussion preceding it.

A second probabilistic interpretation of the density $\Lambda_\balpha$
can be derived from random matrix theory. Consider the diagonal matrix
$D = \diag( a_1, \ldots, a_m )$ with positive non-zero entries
and its orbit $\Omega$ under unitary conjugation in the space of
$m \times m$ positive semidefinite matrices.
If $\mu$ is the normalized $U(m)$-invariant measure on $\Omega$ then 
$\Lambda_\balpha( x ) \intd x$ is the distribution of any diagonal entry of a random
positive semidefinite matrix of arbitrary size distributed according to $\mu$.
See Section~3 and \cite[Section~8]{OV}, or \cite{Faraut}, for details.

The third occurrence of P\'olya frequency functions as characteristic functions
of probability distributions arises indirectly from multivariate statistics
via orbital integrals. Such integrals over matrix groups go back to Wishart's
original work, widely considered to be the origin of modern random
matrix theory. The precise calculation of orbital integrals of the type we
discussed above produced closed-form expressions for various
probability distributions of matricial variables. Of particular interest
is the setting where one collects $n$~independent observations
$X^{(1)}, \ldots, X^{(n)}$ from a $p$-dimensional Gaussian vector
$X \sim N( \mu, \Sigma )$ with mean $\mu \in \R^p$ and covariance
matrix $\Sigma \in \R^{p \times p}$. A fundamental problem of great
interest to applied scientists is to detect non-spurious correlations
between the components of~$X$. If $\bX$ denotes the $n \times p$
matrix whose rows are the vectors $X^{(1)}, \ldots, X^{(n)}$,
then detecting dependencies between the variables can be
performed by computing the associated \emph{sample covariance matrix},
which, up to rescaling, equals $\bX^T \bX$. A large entry in
$\bX^T \bX$ indicates a high level of dependence between the corresponding
variables and understanding the exact distribution of $\bX^T \bX$ is
critical to assessing whether a large entry arises purely by chance or as
the result of a real interaction. Wishart was able to compute the density of
such matrices in closed form via changes of variables and the calculation of Jacobians.
The resulting distribution is named after him. When the observations $X^{(1)}$, \ldots, $X^{(n)}$
have different means $\mu_1$, \ldots, $\mu_n$
(the \emph{non-central} case), computing the density of $\bX^T \bX$
involves integrals over the orthogonal group and the result can be written
in terms of zonal polynomials~\cite{Takemura}. The challenges encountered
by statisticians along this path are well recorded in the monograph by
Farrell~\cite{Farrell}. It seems Karlin himself lectured around 1960 on this
subject, providing a versatile Haar measure disintegration formula.
No doubt elements of \cite{Karlin} were on his desk at that time.

\section{Two atoms: the Berezin--Gindikin--Wallach phenomenon in total positivity}

As seen in Theorem~\ref{TPFprob}, the `atoms' that make up a
general P\'olya frequency function are exponential random variables,
together with at most one Gaussian. If we consider a single atom
then it is evident that any positive real
power of the associated P\'olya frequency function is also in
the same class. However, as we now explain, the picture is markedly
different for more than one exponential variable.

The question of which powers of the density of the sum of two standard
exponential variables are totally non-negative is partially answered in 
the 1964 paper~\cite{KarlinTAMS} of Karlin. His results show that all integer powers of the
associated $\PF$ function $\Omega$ are again $\PF$ functions, but the non-integer powers
are only shown to be $\TN$ of some finite order.

\begin{definition}
Given totally ordered sets $X$ and $Y$, and an integer $p \geq 1$, a kernel
$K : X \times Y \to \R$ is said to be \textit{totally non-negative of order $p$},
denoted $\TN_p$, if for all $x_1 < \cdots < x_r$ in $X$ and $y_1 < \cdots < y_r$ in $Y$,
where $1 \leq r \leq p$, the determinant
$\det (K(x_i, y_j))_{i, j = 1}^r \geq 0$.

A function $\Lambda : \R \to \R$ is said to be $\TN_p$ if its associated
Toeplitz kernel $T_\Lambda$ given by~(\ref{ETlambda}) is $\TN_p$. We say
$\Lambda$ is \textit{totally non-negative} if $\Lambda$ is $\TN_p$ for
all $p \geq 1$.
\end{definition}

With this notation at hand, the result of Karlin can now be stated.
We denote the set of non-negative integers $\{ 0, 1, 2, \ldots \}$ by $\Z_{\geq 0}$.

\begin{theorem}[Karlin, \cite{KarlinTAMS}]\label{Tkarlin}
Let $\Omega : \R \to \R; \ x \mapsto \bone_{x \geq 0} \, x e^{-x}$
be the probability density function for the sum of two independent
standard exponential random variables.
Given an integer $p \geq 2$ and a scalar $\alpha \geq 0$,
the function $\Omega^\alpha : x \mapsto \Omega( x )^\alpha$
is $\TN_p$ if $\alpha \in \Z_{\geq 0} \cup ( p - 2, \infty )$.
\end{theorem}

In particular, if $\alpha$ is a non-negative integer power, then
$\Omega^\alpha$ is $\TN$. This case was shown previously by Schoenberg
in \cite{Schoenberg51} as an immediate consequence of
Theorem~\ref{Tlaplace}.

It is natural to ask if the converse to this 1964 result of Karlin holds.
To the best of our knowledge, this question was not answered in the
literature, and it is the first result that we announce in this note.
In fact, we prove a twofold strengthening.

\begin{theorem}[\cite{Khare2020}]\label{TstrongKarlin}
Given $q$, $r \in (0,\infty)$, let $\Omega_{(q,r)}$ be the
probability density function for $q X_1 + r X_2$,
where $X_1$ and $X_2$ are independent standard exponential
random variables. Now fix an integer $p \geq 2$ and a real number $\alpha \geq 0$.
\begin{enumerate}
\item The function
\[
\Omega_{(q,r)}^\alpha : \R \to \R; \ x \mapsto \Omega_{(q,r)}( x )^\alpha
\]
is $\TN_p$ if $\alpha \in \Z_{\geq 0} \cup ( p - 2, \infty )$.

\item If $\alpha \in ( 0, p - 2 ) \setminus \Z$, then
$\Omega_{(q,r)}^\alpha$ is not $\TN_p$. More strongly, given arbitrary
real numbers $x_1 < \cdots < x_p$ and $y_1 < \cdots < y_p$, there exists
$a \in \R$ such that the matrix
\[
(\Omega_{(q,r)}(x_j - y_k - a)^\alpha)_{j,k=1}^p
\]
is $\TP$ if $\alpha > p-2$, $\TN$ if $\alpha \in \{ 0, 1, \ldots, p-2 \}$,
and has a negative principal minor if $\alpha \in (0,p-2) \setminus \Z$.
\end{enumerate}
\end{theorem}

The proof of this result is obtained by exploiting a variant of
Descartes' rule of signs that was very recently shown and used by
Jain~\cite{Jain2} in her study of entrywise powers preserving positive
semidefiniteness. Note that Theorem~\ref{TstrongKarlin} strengthens
Karlin's result:
(a)~it extends Theorem~\ref{Tkarlin} from $q=r=1$ to
all $q$, $r>0$;
(b)~it extends the total non-negativity in Theorem~\ref{Tkarlin} to total
positivity on a large collection of $p \times p$ matrices, and
(c)~it shows that the converse to the extended
Theorem~\ref{Tkarlin} holds for all $q$, $r>0$.

This last point means that the set of powers preserving the $\TN_p$
property for one-sided P\'olya frequency functions built out of 
two exponential variables is comprised of an arithmetic
progression and a semi-infinite axis. In the subsequent decade to
Karlin's 1964 work, this phenomenon was observed in numerous different
contexts in mathematics.
\begin{itemize}
\item In complex analysis, Rossi and Vergne~\cite{RossiVergne} classified
the powers $\alpha$ of a Bergman kernel on a fixed tube domain
$D \subset \C^n$ which are the reproducing kernels for some Hilbert space
of holomorphic functions on $D$. They termed this set of powers the
\textit{Wallach set}, and showed that it equals an arithmetic progression
together with a half-line, precisely as above.

\item Rossi--Vergne named the above set following Wallach, who was
then following Harish-Chandra and studying holomorphic discrete series of
connected simply connected Lie groups $G$. Wallach classified
in~\cite{Wallach} the twist parameters $\alpha$ of the center of the
maximal compact reductive subgroup $K$ of $G$ for which the
$K$-finite highest weight module over the Lie algebra of $G$ has certain
unitarizability properties. Once again, $\alpha$ belongs to a similar
set.

\item In his pioneering work on quantization, Berezin encountered a
similar set while classifying admissible values of Planck's constant $h$
in deformations of bounded symmetric domains~\cite{Berezin}.
Specifically, \cite[Theorem~1.1]{Berezin} provides a complete picture of
such admissible $h$, adapted to Cartan's four classes of classical
symmetric domains.
As might be expected, these tables display a mixture of continuous and
discrete values for~$h$. Moreover, \cite[Lemma~1.1]{Berezin} is perfectly
aligned with the classification of fractional powers which preserve the
positivity of a structured kernel: \emph{Given a bounded homogeneous
domain with associated Bergman kernel~$K$, for a real number $h$ to be
admissible it is necessary and sufficient that the kernel $K^{1/h}$ is
positive semidefinite.}

\item Gindikin worked in~\cite{Gindikin} with Riesz distributions
$R_\mu$ associated to symmetric cones and indexed by a complex parameter $\mu$
and showed that $R_\mu$ is a positive measure if and only if $\mu$ lies in some similar
set.

\item Finally, FitzGerald and Horn~\cite{FitzHorn} classified the set of
entrywise powers preserving positive semidefiniteness on $p \times p$
matrices, and showed that this set also equals $\Z_{\geq 0} \cup
(p-2,\infty)$.
\end{itemize}

It is remarkable that all of the above named authors obtained their
results in such diverse areas of mathematics within a few years of each other
(all in the 1970s). Each of these works has been followed by
tremendous activity. More recently, such a set has been found
in the theory of non-central Wishart distributions
(see, for example, \cite{FK-book,LetacMassam,Mayer,PR}). Karlin's
work predated all of these works and results, thus providing an earlier
instance of such a Berezin--Gindikin--Wallach set.
We also mention the very recent result of Sra (see \cite[Theorem~2]{Sra}) obtaining such a spectrum from the characterization of fractional powers that preserve the positivity of Hua--Bellman matrices.

In his comprehensive 1968 book~\cite{Karlin} on total
positivity, Karlin provided a second such set of powers.
Let the set of positive integers $\{ 1, 2, 3, \ldots \}$ be denoted by
$\Z_{> 0}$.

\begin{theorem}[Karlin, \cite{Karlin}]\label{Tkarlin2}
If $\Lambda$ is a one-sided $\PF$ function 
and $p \geq 2$ is an integer then $\cB \{ \Lambda \}^\alpha$
is the Laplace transform of a $\TN_p$ function
for all $\alpha \in \Z_{> 0} \cup ( p - 1, \infty )$.
\end{theorem}

Once again, it is natural to ask if the converse holds. We answer this
question in the affirmative, and add a third equivalent condition
involving a single test function.

\begin{proposition}[\cite{Khare2020}]
Fix an integer $p \geq 2$ and a scalar $\alpha \geq 1$. The following are
equivalent:
\begin{enumerate}
\item If $\Lambda$ is a one-sided $\PF$ function, then
$\cB \{ \Lambda \}^\alpha$ is the Laplace transform of a $\TN_p$ function.

\item For the `single atom' $\PF$ function
\[
\lambda_1  : \R \to \R; \ x \mapsto  \bone_{x \geq 0} \, e^{-x},
\]
the density of a standard exponential variable, the power
$x \mapsto \cB \{ \lambda_1 \}( x )^\alpha$ is the Laplace transform of a $\TN_p$
function.

\item $\alpha \in \Z_{>0} \cup ( p - 1, \infty )$.
\end{enumerate}
\end{proposition}

That $(3) \implies (1)$ was Karlin's 1968 result above, and
is a consequence of Schoenberg's representation Theorem~\ref{Tlaplace}
for $\PF$ functions. (The implications $(1) \implies (2) \implies (3)$
are similarly not hard to prove.)

We close this section with a few results related to the topics addressed above.
The first is from Schoenberg's 1955 work~\cite{Schoenberg55} that initiated
the study of $\TN_p$ kernels (which Schoenberg termed
\textit{multiply positive functions}). In this work, Schoenberg studied the order of total
non-negativity of powers of the \textit{Wallis kernel}.

\begin{theorem}[Schoenberg, \cite{Schoenberg55}]
Let
\[
W : \R \to \R; \ x \mapsto \bone_{| x | \leq \pi / 2} \, \cos x.
\]
For any integer $p \geq 2$, the power $W^\alpha$ is $\TN_p$
if and only if $\alpha \geq p - 2$.
\end{theorem}

In analogy with Theorem~\ref{TstrongKarlin}, we strengthen the total
non-negativity to total positivity, and the lack thereof to principal
minors, on large subsets of arguments.

\begin{theorem}[\cite{Khare2020}]
Let $p \geq 2$ be an integer.
Given arbitrary reals $x_1 < \cdots < x_p$ and $y_1 < \cdots < y_p$,
there exists a `multiplicative shift' $m_0 \in (0,\infty)$ such that the
matrices
\[
(W(m(x_j - y_k))^\alpha)_{j,k=1}^p, \qquad 0 < m < m_0
\]
are each $\TP$ if $\alpha > p-2$, $\TN$ if $\alpha \in \{ 0, 1, \ldots, p - 2 \}$,
and each have a negative principal minor if $\alpha \in ( 0, p - 2 ) \setminus \Z$.
\end{theorem}

The other result which we address here is the classification of the
discontinuous P\'olya frequency functions. In
his 1951 paper~\cite{Schoenberg51}, as well as the preceding
announcements in \textit{Proc.\ Natl.\ Acad.\ Sci.\ USA}, Schoenberg
asserts that the only discontinuous P\'olya frequency function is the
standard exponential density $\lambda_1( x ) := \bone_{x \geq 0} \, e^{-x}$
``up to changes in scale and origin'' (which includes reflecting about the $y$~axis).
This also implies that the only discontinuous totally non-negative
function is the Heaviside function $H_1( x ) := \bone_{x \geq 0}$, up
to changes in scale and origin and multiplying by an exponential
factor $e^{a x + b}$, where $a$, $b \in \R$.

It turns out that these statements are not quite true, precisely at the
point of discontinuity. Our explorations in~\cite{BGKP-TN}
led us to a family of discontinuous $\PF$ functions $\{ \lambda_d : d \in [0,1] \}$
that lies outside the above class, and subsequently, to a small
correction of Schoenberg's classification.

\begin{theorem}[\cite{Khare2020}]\label{Tdiscont}
A P\'olya frequency function is discontinuous if and only if it equals,
up to changes in scale and origin,
the following function $\lambda_d$ for some $d \in [ 0, 1 ]$:
\[
\lambda_d : \R \to \R; \ x \mapsto \begin{cases}
e^{-x} \qquad & \text{if } x > 0,\\
d \qquad & \text{if } x = 0,\\
0 \qquad & \text{if } x < 0.
\end{cases}
\]
Similarly, a $\TN$ function is discontinuous if and only if (up to
changes in scale and origin) it equals $x \mapsto e^{a x + b} \lambda_d( x )$
for some $d \in [0,1]$ and $a$, $b \in \R$.
\end{theorem}

Note that the Laplace transform of $\lambda_d$ does not depend on~$d$.

\section{Three or more atoms: non-smooth P\'olya frequency functions as
hypoexponential densities}

We now turn from Berezin--Gindikin--Wallach type sets -- encountered here
in the study of $\TN_p$ powers of functions  --  to transforms of
broader classes of P\'olya frequency
functions. The preceding result classifies the discontinuous $\PF$
functions, following Schoenberg~\cite{Schoenberg51}. In the same paper,
Schoenberg classified the non-smooth P\'olya frequency functions.
In the language of probability theory, they are precisely the
\textit{densities of finite linear combinations of independent exponential variables}.
Such functions were studied in detail by Hirschman and Widder, first in their
1949 paper~\cite{HW49}, and then in the 1955 memoir~\cite{HW}. In recent
work~\cite{BGKP-HW}, we further investigate these maps, obtaining
connections to probability and to symmetric function theory; the goal of
this section is to announce some of those results. Here and in ~\cite{BGKP-HW}, we term these frequency functions
\textit{Hirschman--Widder densities}.

Before we discuss the power preservers of such densities (following
Karlin, Schoenberg, and the results in the preceding section), we present
some pleasing properties of Hirschman--Widder densities. The
preceding section dealt with single exponential random variables which
led to discontinuous $\PF$ functions; thus, in this
section we consider the densities of variables
\begin{equation}\label{Eerlang}
\alpha_1 X_1 + \cdots + \alpha_m X_m, \qquad \text{where } m \geq 2
\text{ and } \balpha := ( \alpha_1, \ldots, \alpha_m ) \in ( \R^\times )^m.
\end{equation}
Here $X_1$, \ldots, $X_m$ are independent standard exponential random
variables and $\R^\times$ denotes the set of non-zero real numbers.

We denote the corresponding density function by
$\Lambda_\balpha$. Such functions are known in
probability and statistics as \textit{hypoexponential densities}, 
or as \textit{Erlang densities} if the coefficients $\alpha_j$ are all equal,
and are relevant to several applied fields. However, the
connection to the work of Hirschman and Widder seems to not be widely
known in the probability literature.

Here are some of the properties enjoyed by Hirschman--Widder densities.
\begin{enumerate}
\item The function $\Lambda_\balpha : \R \to [0,\infty)$ is
the unique continuous function with bilateral Laplace transform
\begin{equation}
\cB \{ \Lambda_\balpha \}(s) = \prod_{j=1}^m \frac{1}{1 + \alpha_j s}
\end{equation}
on the open half-plane
$\{ s \in \C : \Re s > -\alpha_j^{-1} \text{ for } j = 1, \ldots, m \}$.

\item In particular, $\Lambda_\balpha$ is both the probability density
function of the random variable (\ref{Eerlang}) and a P\'olya frequency
function.

\item $\Lambda_\balpha$ has an additive representation
via one-sided exponential densities, as well as a multiplicative one via
the bilateral Laplace transform, which corresponds to the
convolution of these densities.
\end{enumerate}

The additive representation is particularly gratifying when the parameters
$\alpha_j$ are pairwise distinct and positive:
\begin{equation}
\Lambda_\balpha(x) = \bone_{x \geq 0} \sum_{j=1}^m a_j
e^{-a_j x} \prod_{k \neq j} \frac{a_k}{a_k - a_j}, 
\qquad \text{where } a_j := \alpha_j^{-1} \text{ for all } j.
\end{equation}

\subsection{Taylor coefficients, moments, and symmetric
functions}\label{SHW}

We now turn to some connections between
Hirschman--Widder densities and symmetric function theory. As discussed
above, Schoenberg's results in~\cite{Schoenberg51}
imply that the densities $\Lambda_\balpha$ are the
\textit{only} one-sided, non-smooth, continuous P\'olya frequency
functions that vanish on $( -\infty, 0 )$ and are non-zero on $( 0, \infty )$.

It is natural to take a closer look at the Taylor expansion of
$\Lambda_\balpha$ at $0^+$ and also at the moments of this density
function. The descriptions of both of these quantities
involve a well-known family of symmetric functions, the
\textit{complete homogeneous symmetric polynomials}:
\begin{equation}
h_p( a_1, \ldots, a_m ) := \sum_{1 \leq j_1 \leq j_2 \leq \cdots \leq j_p
\leq m} a_{j_1} a_{j_2} \cdots a_{j_p}, \qquad \text{for all } p \in
\Z_{> 0}.
\end{equation}
(We set $h_0 \equiv 1$.)
These polynomials are clearly symmetric in their arguments, and equal
Schur polynomials corresponding to tableaux for a single row with $p$
cells and an alphabet of size $m$. In particular, they carry
representation-theoretic content: they correspond to the characters of
certain polynomial representations of the symmetric group, or of the
special Lie algebra $\mathfrak{sl}_m(\C)$.

With these polynomials at hand, we provide closed-form expressions for
the Maclaurin coefficients at $0^+$, as well as the moments, of
Hirschman--Widder densities.

\begin{theorem}[\cite{BGKP-HW}]
Fix an integer $m \geq 2$ and parameters
$\balpha = ( \alpha_1, \ldots, \alpha_m ) \in ( 0, \infty )^m$.
\begin{enumerate}
\item The Hirschman--Widder density $\Lambda_\balpha$
has the expansion
\begin{equation}\label{Emaclaurinhw}
\Lambda_\balpha( x ) = \alpha_1 \cdots \alpha_m %
\sum_{n = m - 1}^\infty %
\frac{( -1 )^{n - m + 1} h_{n - m + 1}( \alpha_1^{-1}, \ldots,
\alpha_m^{-1} )}{n!} x^n
\end{equation}
as its Maclaurin series, convergent for all $x \in [0,\infty)$. In
particular, $\Lambda_\balpha$ is smooth except at the origin, where it is
of continuity class $C^{m-2}$ but not $C^{m-1}$.

\item The Hirschman--Widder density $\Lambda_\balpha$
has $p^{th}$ moment
\[
\mu_p := \int_\R x^p \Lambda_\balpha(x) \intd x =
p! \, h_p(\alpha_1, \ldots, \alpha_m)  \qquad \text{for all } p \geq 0.
\]
\end{enumerate}
\end{theorem}

While both the Maclaurin coefficients and the moments are
expressed by the same family of symmetric functions, the arguments are
different, involving $\{ a_j = \alpha_j^{-1} \}$ and $\{ \alpha_j \}$,
respectively. Corresponding
to these are two familiar identities which appear as
byproducts of the proofs in~\cite{BGKP-HW}. The
first is the `local' version of the
well-known generating function for the symmetric polynomials $h_p$:
\[
\sum_{p = 0}^\infty h_p( \alpha_1^{-1}, \ldots, \alpha_m^{-1} ) z^p = %
\prod_{j = 1}^m \frac{1}{1 - \alpha_j^{-1} z} \quad \text{whenever } %
| z | < \min\{ \alpha_j : j = 1, \ldots, m \}.
\]
The second is the moment-generating function for the density
$\Lambda_\balpha$ of $X := \sum_{j=1}^m \alpha_j X_j$:
\[
\sum_{p = 0}^\infty \frac{\mu_p}{p!} z^p = \E[ e^{z X} ] = %
\cB\{ \Lambda_\balpha \}( -z ) = %
\prod_{j = 1}^m \frac{1}{1 - \alpha_j z}.
\]

As a final result in this vein, we show in~\cite{BGKP-HW} that for each
$m \geq 2$, the parameters $\balpha$ can be recovered from finitely many
moments or Maclaurin coefficients.

\begin{theorem}[\cite{BGKP-HW}]
Given $m \geq 2$ and $\balpha \in (0,\infty)^m$, the
parameters $\balpha$ may be recovered, up to permuting its
entries, from the first $m$ moments $\mu_1$, \ldots, $\mu_m$ of the
density $\Lambda_\balpha$, and also from the lowest $m + 1$
non-trivial Maclaurin coefficients of $\Lambda_\balpha$, that is,
$\Lambda_\balpha^{(k)}(0^+)$ with $m - 1 \leq k \leq 2m - 1$.
\end{theorem}

The proof crucially relies upon the Jacobi--Trudi identity from symmetric
function theory. As we mention below, this is one of several hitherto
unexplored connections between symmetric functions and positivity that
have emerged in recent works.

\subsection{Power and polynomial preservers of Hirschman--Widder
densities}

We now return to the theme of understanding
which transformations preserve the class of
Hirschman--Widder densities. As Karlin and Schoenberg have shown, the
density of $X_1 + X_2$ (with $X_1$ and $X_2$ independent standard
exponential random variables) is a $\PF$ function, every positive
integer power of which is also a $\PF$ function.
As we asserted in Theorem~\ref{TstrongKarlin}, the situation is identical
for $\alpha_1 X_1 + \alpha_2 X_2$ whenever $\alpha_1$, $\alpha_2 > 0$. More
generally, we have the following result;
note that, since $\Lambda_\balpha$ is symmetric in
$\balpha$, the entries of $\balpha$ may be assumed to be monotonic.

\begin{theorem}[\cite{BGKP-HW}]
Consider the class of probability densities
\[
\{ \Lambda_\balpha : \alpha_1^{-1}, \ldots, \alpha_m^{-1}
\text{ are positive and form an arithmetic progression} \}
\]
of finite linear combinations $\sum_{j=1}^m \alpha_j X_j$ of independent
standard exponential random variables. This class of densities is closed under
taking positive integer powers.
\end{theorem}

In particular, if $m=2$, then all positive integer powers of
$\Lambda_\balpha$ are Hirschman--Widder densities, so P\'olya frequency functions,
so totally non-negative functions.

However, the preceding result is somewhat misleading as far as general
parameters~$\balpha$ go. For $m \geq 3$ and almost all
$\balpha$, we show that not only do integer powers not remain
$\PF$ functions, or even $\TN$ functions, but in fact \textit{no}
polynomial transform enjoys this property, save for the obvious
ones.

\begin{theorem}
Fix an integer $m \geq 3$. There exists a set
$\cN \subset (0,\infty)^m$ with zero Lebesgue measure
such that
\[
p \circ \Lambda_\balpha : \R \to \R; \ x \mapsto p( \Lambda_\balpha( x ) )
\]
is not a $\TN$ function, and so not a P\'olya frequency function,
for any $\balpha \in ( 0, \infty )^m \setminus \cN$ and any
non-homothetic real polynomial $p$, that is, any real polynomial not
of the form $p( x ) = c x$ for some $c > 0$.
\end{theorem}

\section{One and two-sided P\'olya frequency functions and sequences}

Having studied finite linear combinations of exponential random
variables and the transformations that preserve total non-negativity
of their densities, we now tackle the same preserver problem for the
general case: one-sided and two-sided P\'olya frequency functions.
This is the focus of the third paper~\cite{BGKP-TN} that is being
announced in the present note.

In studying these preservers, we were
naturally motivated by the celebrated paper of
Schoenberg~\cite{Schoenberg42}, which classified the preservers of a
related notion: positive semidefiniteness. Coming from considerations of
distance geometry (specifically, positive definite functions of the
form $F \circ \cos$ on Euclidean spheres) Schoenberg
classified the continuous functions $F : [-1,1] \to \R$ which
preserve positivity when applied entrywise to positive semidefinite
matrices of all dimensions.

\begin{theorem}[Schoenberg, \cite{Schoenberg42}]\label{Tspheres}
Suppose $F : [-1,1] \to \R$ is continuous. The following are equivalent.
\begin{enumerate}
\item The \textit{entrywise transform} $F[A] := (F(a_{ij}))_{i,j=1}^n$ is positive
semidefinite for every positive semidefinite matrix  $A = (a_{ij})_{i,j=1}^n$
with entries in $[-1,1]$ and every $n \geq 1$.

\item $F$ admits a power-series representation
$\sum_{k=0}^\infty c_k x^k$ on $[-1,1]$ with non-negative Maclaurin
coefficients.
\end{enumerate}
\end{theorem}

As observed by P\'olya and Szeg\H{o} in 1925, the implication
$(2) \implies (1)$ essentially follows from the Schur product
theorem~\cite{Schur1911}, and functions which such a representation are called \textit{absolutely monotonic}.
Schoenberg's theorem \cite{Schoenberg42}
provides the highly non-trivial converse result, that continuous preservers are precisely the absolutely monotonic functions.
The continuity assumption in Theorem~\ref{Tspheres} was subsequently removed by
Rudin, who showed moreover that the test set can be greatly reduced in every
dimension, to low-rank Toeplitz matrices.

\begin{theorem}[Rudin, \cite{Rudin59}]\label{Trudin}
Suppose $I = (-1,1)$ and $F : I \to \R$. The following are equivalent.
\begin{enumerate}
\item $F[A]$ is positive semidefinite for every positive semidefinite
 $A \in I^{n \times n}$ and every $n \geq 1$.

\item $F[A]$ is positive semidefinite for every positive semidefinite
 $A \in I^{n \times n}$ which is Toeplitz of rank at most~$3$ and every $n \geq 1$.

\item $F$ is represented by a power series $\sum_{k=0}^\infty c_k x^k$
on $I$, with all $c_k \geq 0$.
\end{enumerate}
\end{theorem}

Inspired by Schoenberg's theorem, and with applications to moment
problems in mind, we recently strengthened Theorem~\ref{Tspheres} in a
parallel direction to Rudin's result.

\begin{theorem}[\cite{BGKP-hankel}]\label{Thankel}
The three conditions in Theorem~\ref{Trudin} are also equivalent to
the following.
\begin{enumerate}
\setcounter{enumi}{3}
\item $F[A]$ is positive semidefinite for every positive semidefinite
$A \in I^{n \times n}$ which is Hankel of rank at most~$3$ and every $n \geq 1$.
Equivalently, $F[-]$ preserves the class of
moment sequences of positive measures on $I$ having moments of all orders.
\end{enumerate}
\end{theorem}

These results can be recast in the language of composition maps $C_F$. A
real $n \times n$ matrix $A$ can be identified with a function
$A : [n] \times [n] \to \R$, where $[n] := \{ 1, \ldots, n \}$.
Now, the entrywise transform $F[A]$ is simply the composition
$C_F(A) := F \circ A$. Thus, the above results are equivalent to
classifying the composition maps $C_F$ that preserve the class of
\textit{positive semidefinite kernels} on an infinite domain $X$. In
particular, Theorem~\ref{Thankel} classifies the composition maps that
preserve positivity on the family of real Hankel kernels on
$\Z_{\geq 0}^2$.

It is natural to move from studying the preservers of structured
Hankel kernels to those of structured Toeplitz kernels. In this
section, we consider endomorphisms of several different classes of
functions: P\'olya frequency functions, totally non-negative
functions, and P\'olya frequency sequences.  We also discuss one-sided
variants of these. We begin by asserting that the classes of P\'olya
frequency functions and sequences are as rigid as the class of
Hirschman--Widder densities, not just for polynomials, but for
arbitrary composition transforms.

\begin{theorem}[\cite{BGKP-TN}]\label{TPF}
Let $F : [ 0, \infty ) \to [ 0, \infty )$ be a function. Consider the
composition transform $C_F$, taking
any function of the form
$\Lambda : X \to [ 0, \infty )$ 
to the map $C_F( \Lambda ) := F \circ \Lambda : X \to [ 0, \infty )$.
\begin{enumerate}
\item The composition transform $C_F$ preserves the class of
P\'olya frequency functions if and only if $F$ is a positive
homothety: there exists $c > 0$ such that $F( x ) = c x$ for all $x$.

\item The transform $C_F$ preserves the class of P\'olya
frequency sequences if and only if $F$ is a
positive homothety or a non-negative constant:
either $F \equiv 0$ or there exists $c > 0$ such that $F \equiv c$
or $F( x ) = c x$ for all $x$.

\item The transform $C_F$ preserves the class of totally non-negative
functions if and only if $F \equiv 0$ or
there exists $c > 0$ such that $F \equiv c$ or $F( x ) = c x$ for all
$x$ or $F( x ) = c \bone_{x > 0} $ for all $x$.
\end{enumerate}
\end{theorem}

We next present a one-sided version of Theorem~\ref{TPF}. Given
a domain $X \subset \R$, we say that a function $\Lambda : X \to \R$ is
\textit{one sided} if there exists $x_0 \in \R$ such that
$\Lambda \equiv 0$ on $X \cap ( -\infty, x_0 )$ or
$\Lambda \equiv 0$ on $X \cap ( x_0, \infty )$.

\begin{theorem}[\cite{BGKP-TN}]
Suppose $F : [ 0, \infty ) \to [ 0, \infty )$.
\begin{enumerate}
\item The composition transform $C_F$ preserves the
class of one-sided P\'olya frequency functions if and only if
there exists $c > 0$ such that $F( x ) = c x$ for all $x$.

\item The transform $C_F$ preserves the class of
one-sided P\'olya frequency sequences if and only if
$F \equiv 0$ or there exists $c > 0$ such that
$F( x ) = c x$ for all $x$.

\item The transform $C_F$ preserves the class of
one-sided totally non-negative functions if
and only if $F \equiv 0$ or there exists $c > 0$ such that
either
$F( x ) = c x$ for all $x$ or $F( x ) = c \bone_{x > 0}$ for all $x$.
\end{enumerate}
\end{theorem}

These results reveal the rigidity of the endomorphisms of several
classes of totally non-negative functions. A similar phenomenon is
uncovered when understanding the endomorphisms of totally positive
functions, or of $\TP$ P\'olya frequency functions or
sequences. Recall that these are the subsets of the corresponding
classes of maps for which all determinants in the defining conditions
are positive instead of just non-negative.

\begin{theorem}[\cite{BGKP-TN}]\label{TTPPF}
Suppose $F : (0,\infty) \to (0,\infty)$. The composition
transform $C_F$
preserves any of the following classes (and so all of them)
if and only if there exists $c > 0$ such that
$F( x ) = c x$ for all $x$: (a)~totally positive $\PF$ functions,
(b)~totally positive $\PF$ sequences,
(c)~totally positive functions.
\end{theorem}

The proof of Theorem~\ref{TTPPF} involves two interesting ingredients.
The first is a density phenomenon for totally positive P\'olya
frequency functions.

\begin{proposition}[\cite{BGKP-TN}]
The class of totally positive P\'olya frequency functions is dense in the
class of P\'olya frequency functions $\Lambda$ that are \emph{regular},
that is, such that
$\Lambda( x ) = ( \Lambda( x^+ ) + \Lambda( x^-) ) / 2$ for all
$x$, where $\Lambda( x^+ )$ and $\Lambda( x^- )$ are the right-hand and
left-hand limits of the function $\Lambda$ at the point $x$.
\end{proposition}

The reader may be reminded of the density result of Anne M.\ Whitney
\cite{Whi}, which asserts that the set of $m \times n$ $\TP_p$ matrices
is dense in that of $m \times n$ $\TN_p$ matrices, for all
positive integers $m$, $n$ and $p$. The proposition
above is the analogous result for Toeplitz kernels arising from
P\'olya frequency functions.

The second ingredient in its proof is the discontinuous test function
$\lambda_{1/2}$ from Theorem~\ref{Tdiscont}, which is indeed a regular
P\'olya frequency function. In fact, it was in this context that we
first encountered $\lambda_{1/2}$ and were led to the more general
family of functions $\lambda_d$ and to the classification
Theorem~\ref{Tdiscont}.

\section{Preservers of total positivity on arbitrary domains}

In this final section, we go from working with structured kernels to
arbitrary ones. Just as the preservers of all positive semidefinite
matrices are classified by Schoenberg's theorem~\ref{Tspheres} and its
strengthenings in \cite{Rudin59} and \cite{BGKP-hankel}, a natural and
parallel question is to understand the preservers of total
positivity. Here, we present the complete solution to the following
problems.

\begin{quote}
\textit{Given non-empty totally ordered sets $X$ and $Y$ and a
positive integer~$p$,
classify the composition transforms $C_F$ which preserve the classes of
$\TN_p$, $\TP_p$, $\TN$, and $\TP$ kernels on $X \times Y$.}
\end{quote}

\begin{remark}
Two initial observations are in order.
First, from the definitions, it is clear that we can assume $X$ and
$Y$ have cardinality at least~$p$, for the $\TN_p$ and $\TP_p$
problems, and are infinite otherwise, to ensure the relevant class is
not empty.
Second, the kernel $K \equiv 0$ is $\TN$ and $K' \equiv 1$ is $\TP_1$
on any non-empty domain $X \times Y$, so there are no further 
restrictions for the $\TN_p$ and $\TP_1$ classification problems.
\end{remark}

As in \cite{BGKP-TN}, we examine the totally
non-negative and totally positive cases separately.

\subsection{Preservers of total non-negativity}

It is immediate that the functions $F$ such that $C_F$ preserves the
$\TN_1$ or $\TP_1$ property for kernels are precisely the self maps
$F : [ 0, \infty ) \to [ 0, \infty )$ or
$F : ( 0, \infty ) \to ( 0, \infty )$ respectively. Thus, in the
sequel we will assume $p \geq 2$.

To write down systematically all total non-negativity preservers, we
first introduce the following compact notation.

\begin{definition}
Given totally ordered sets $X$ and $Y$, and a positive integer $p$,
let
\begin{align*}
\sFTN_{X, Y}(p) := &\ \{ F : [ 0, \infty ) \to \R \mid %
F \circ K \text{ is $\TN_p$ for any $\TN_p$ kernel $K$ on $X \times Y$} \}\\
\sFTP_{X, Y}(p) := &\ \{ F : [ 0, \infty ) \to \R \mid %
F \circ K \text{ is $\TP_p$ for any $\TP_p$ kernel $K$ on $X \times Y$}
\}.
\end{align*}
We also let $\sFTN_{X, Y}(\infty)$ and $\sFTP_{X, Y}(\infty)$ denote
the analogous classes of functions for $\TN$ and $\TP$ kernels,
respectively, and use $\TN_\infty$ and $\TP_\infty$ synonymously with
$\TN$ and $\TP$.
\end{definition}

To classify the total non-negativity preservers, we first note that if
$X$ and $Y$ have cardinality at least~$p$ then every preserver of the
set of $\TN_p$ matrices in $\R^{p \times p}$ is automatically a
preserver of the set of $\TN_p$ kernels on $X \times Y$, and vice
versa. Indeed, one inclusion of preservers is immediate, and the
reverse inclusion follows by padding a $p \times p$ $\TN_p$ matrix by
zeros to yield a $\TN_p$ kernel on $X \times Y$. With this observation
at hand, we present the solution to the easier of the two questions
posed above.

\begin{theorem}[\cite{BGKP-TN}]\label{Ttnmatrix}
Let $X$ and $Y$ be totally ordered sets with
cardinality at least $p$, where $2 \leq p \leq \infty$. Then
\begin{enumerate}
\item $\sFTN_{X, Y}(2) = \{ c x^\alpha : c > 0, \ \alpha > 0 \} \cup
\{ c : c \geq 0 \} \cup \{ c \bone_{x > 0} : c > 0 \}$,
\item $\sFTN_{X, Y}(3) = \{ c x^\alpha : c > 0, \ \alpha \geq 1 \} \cup
\{ c : c \geq 0 \}$, and
\item $\sFTN_{X, Y}(p) = \{ c x : c > 0 \} \cup \{ c : c \geq 0 \}$ if
$4 \leq p \leq \infty$.
\end{enumerate}
\end{theorem}

We briefly remark that proving Theorem~\ref{Ttnmatrix} involves
showing that every preserver of $\TN_p$ preserves $\TN_2$, whence is
measurable. Now the problem reduces to solving the multiplicative
Cauchy functional equation in this setting, to deduce that the
preserver is, up to rescaling, a power function or $\bone_{x > 0}$.

\begin{remark}
Interestingly, the preserver $c \bone_{x > 0}$ also features in
Theorem~\ref{TPF}, and the preservers of $\PF$ sequences, which may be
thought of as bi-infinite $\TN$ Toeplitz matrices, coincide with those of
$\TN_p$ kernels for $4 \leq p \leq \infty$.
\end{remark}

\begin{remark}
As Theorem~\ref{Ttnmatrix} shows, there are very few preservers of
total non-negativity of order at least $4$. The situation was the same
when we considered the preservers of $\TN$ Toeplitz kernels on
infinite domains such as $\R$ or $\Z$. However, this changes
dramatically upon going from the Toeplitz to the Hankel setting: there are
large classes of preservers of Hankel $\TN$ matrices of a fixed
dimension or even of all dimensions. We illustrate this with a few
concrete examples.

\begin{enumerate}
\item The map $F : [ 0, \infty ) \to [ 0, \infty )$
is such that the composition transform $C_F$ preserves all $\TN$
Hankel matrices of all dimensions precisely when it equals a
convergent power series $\sum_{k = 0}^\infty c_k x^k$ on
$( 0, \infty )$ with all $c_k \geq 0$ and is such that
$0 \leq F(0) \leq F(0^+)$. This is similar to the class of preservers
in Schoenberg's Theorem~\ref{Tspheres}, and far larger than the
collection of constants and homotheties in Theorem~\ref{Ttnmatrix}.

\item The set of entrywise powers that preserve the class of Hankel
$\TN_p$ matrices coincides~\cite{BGKP-hankel} with the
Berezin--Gindikin--Wallach-type set $\Z_{\geq 0} \cup ( p - 2, \infty )$
obtained by Karlin: see Theorems~\ref{Tkarlin} and~\ref{TstrongKarlin}.
The same set of powers also appears as the set of power preservers of
the class of continuous Hankel $\TN_p$ kernels on arbitrary
sub-intervals of $\R$ with positive measure, as shown
in~\cite{Khare2020}. Once again, this set differs vastly from the one
in Theorem~\ref{Ttnmatrix}.

\item In a fixed dimension $p \geq 2$, there are polynomial
preservers with negative coefficients, as obtained
in~\cite{BGKP-fixeddim, KT}. This line of
investigation uncovered unexpected connections between positivity and
Schur polynomials, which preceded the findings in Section~\ref{SHW}:
see~\cite{BGKP-HW,Khare2022}.
\end{enumerate}
\end{remark}

\subsection{Preservers of finite-order total positivity}

We now come to the preservers of totally positive kernels. There are
two distinct cases, when $p$ is finite and when $p$ is infinite, and
each has its own subtleties.

The first observation is that if $X$ and $Y$ each have cardinality at
least $2$ and there exists a $\TP_2$ kernel on $X \times Y$, then $X$
and $Y$ necessarily embed into the real line.

\begin{proposition}[\cite{BGKP-TN}]\label{Pembed}
Let $X$ and $Y$ be non-empty totally ordered sets.
The following are equivalent:
\begin{enumerate}
\item There exists a totally positive kernel on $X \times Y$.

\item There exists a $\TP_2$ kernel on $X \times Y$.

\item At least one of $X$ or $Y$ is a singleton set, or there exist
order-preserving embeddings from $X$ into $( 0, \infty )$ and from
$Y$ into $( 0, \infty )$.
\end{enumerate}
\end{proposition}

Thus, the general problem of classifying total positivity preservers
reduces to considering domains embedded in the real line. A key
difference with the corresponding problem for total non-negativity
preservers is that one can no longer extend kernels on smaller domains
to the whole of $X \times Y$ by padding with zeros. Thus, $\TP_p$
preservers of $p \times p$ matrices will preserve $\TP_p$ kernels on
$X \times Y$ whenever $X$ and $Y$ have cardinality at least $p$, but
to show the converse we need to develop additional tools.

The following result is the complete classification in the
finite-order case.

\begin{theorem}[\cite{BGKP-TN}]\label{Ttpmatrix}
Suppose $X$ and $Y$ are totally ordered sets
such that there exists a $\TP_p$ kernel on $X \times Y$ for some
integer $p \geq 2$. Then
\begin{enumerate}
\item $\sFTP_{X, Y}(2) = \{ c x^\alpha : c > 0, \ \alpha > 0 \}$,

\item $\sFTP_{X, Y}(3) = \{ c x^\alpha : c > 0, \ \alpha \geq 1 \}$, and

\item $\sFTP_{X, Y}(p) = \{ c x : c > 0 \}$ if $p \geq 4$.
\end{enumerate}
\end{theorem}

This result should be compared with Theorem~\ref{Ttnmatrix}. The proof
involves first showing that every preserver is continuous; we elaborate
on this below. Once this is done, when $X$ and $Y$ are finite, we may
invoke Whitney's density theorem~\cite{Whi}, which reduces the
classification of preservers of $\TP_p$ kernels on $X \times Y$
to the case of $\TN_p$ kernels, which was done in
Theorem~\ref{Ttnmatrix}. Now some matrix analysis completes the proof.

For finite $X$, $Y$ and $p$, we still need to show that every
preserver is continuous. In fact, a stronger result holds.

\begin{proposition}[\cite{BGKP-TN}]\label{Pvdm}
Let $X$ and $Y$ be totally ordered sets such that there
exists a $\TP_2$ kernel on $X \times Y$.
\begin{enumerate}
\item Every $2 \times 2$ $\TP$ matrix can be embedded in arbitrary
position into a $\TP$ kernel on $X \times Y$.
\item If $F : (0,\infty) \to (0,\infty)$ is
such that the composition map $C_F$ preserves the class of $\TP_p$
kernels on $X \times Y$ for some $p \geq 2$ (including $p = \infty$)
then $F$ is continuous and there exist $c > 0$ and $\alpha > 0$ such
that $F( x ) = c x^\alpha$ for all $x$.
\end{enumerate}
\end{proposition}

The proof of part~(1) involves showing that every $2 \times 2$ $\TP$
matrix is, up to scaling and taking the transpose, a generalized
Vandermonde matrix, and then invoking Proposition~\ref{Pembed}(3).
Moreover, a generalized Vandermonde matrix can easily be embedded in
the $\TP$ kernel
\[
\R \times \R \to \R; \ (x , y) \mapsto \exp(x y),
\]
from which part~(2) follows easily.

This completes the outline of the proof of Theorem~\ref{Ttpmatrix} when $X$
and $Y$ are finite.  However, the proof techniques used for the
results above no longer suffice if at least one of $X$ and~$Y$ is
infinite, even after we know that $F( x ) = c x^\alpha$ as above. The
missing ingredient that is required is an analogue of Whitney's
density theorem for infinite domains. Thus, we extend that result.

\begin{theorem}[\cite{BGKP-TN}]
Let $X$ and $Y$ be subsets of the real line which each contain at
least two distinct points. If $K : X \times Y \to \R$ is a bounded
$\TN_p$ kernel for some integer $p \geq 2$ then there exists a
sequence of $\TP_p$ kernels on $X \times Y$ which converges  to $K$
locally uniformly on the set of points in $X \times Y$ where $K$ is
continuous.
\end{theorem}

\subsection{Total-positivity preservers on bi-infinite domains}

If $X$ or $Y$ is finite, then the $\TP$ kernels on $X \times Y$ are
precisely the $\TP_p$ kernels, where $p$ is the minimum of the
cardinalities of $X$ and $Y$. The preservers in this setting are
classified by Theorem~\ref{Ttpmatrix}. Thus, the only remaining case
requires $X$, $Y$ and $p$ to be infinite.

\begin{theorem}[\cite{BGKP-TN}]\label{Tgeneral}
If $X$ and $Y$ are infinite totally ordered sets such that there
exists a $\TP_2$ kernel on $X \times Y$ then the only preservers of
$\TP$ kernels on $X \times Y$ are the positive homotheties: functions
of the form $F( x ) = c x$ for some $c > 0$.
\end{theorem}

In conjunction with Theorem~\ref{Ttpmatrix}, this result completely
resolves the problem of classifying total positivity preservers for
kernels on an arbitrary domain $X \times Y$.

We briefly sketch here the arguments in the proof. The idea is to use
the classification of the preservers of P\'olya frequency functions
and sequences from a previous section. To do so, one requires a novel
ingredient: order-preserving embeddings of sets containing arbitrarily
long arithmetic progressions. Thus, we introduce the following
definition.

\begin{definition}
Two sets $X$ and $Y$ of real numbers are said to form
\textit{an admissible pair} if for each integer $n \geq 2$ there exist
$n$-step arithmetic progressions
\[
x_1 < \cdots < x_n \text{ in } X \quad \text{and} \quad
y_1 < \cdots < y_n \text{ in } Y
\]
that are equally spaced: $x_{j + 1} - x_j = y_{j + 1} - y_j$ for
$j = 1, \ldots, n - 1$.

Also, let the Minkowski difference
$X - Y := \{ x - y : x \in X, \ y \in Y \}$,
and call a kernel $K : X \times Y \to \R$ \textit{Toeplitz} if there
exists a function $f : X - Y \to \R$ such that
$K( x, y ) = f( x - y )$ for all $x \in X$ and $y \in Y$.
\end{definition}

Now the first step in proving Theorem~\ref{Tgeneral} is the following
common extension of parts (a) and (b) of Theorem~\ref{TTPPF}.

\begin{proposition}[\cite{BGKP-TN}]\label{Padmissible}
Suppose $X$, $Y \subset \R$ form an admissible pair, and let
$F : ( 0, \infty) \to ( 0, \infty )$. The composition map $C_F$ preserves
total positivity for all Toeplitz kernels on $X \times Y$ if and only if
there exists $c > 0$ such that $F( x ) = c x$ for all $x$.
\end{proposition}

With this result at hand, the proof of Theorem~\ref{Tgeneral} goes as
follows. The set $X$ contains an infinite ascending chain or an infinite
descending chain. We construct an order-preserving bijection
$\varphi_X : X \to X'$, where $X' \subset \R$ is a subset that contains
an arithmetic progression of length $2^n$ and step-size $4^{-n}$ for each
$n \geq 2$, and similarly for $\varphi_Y : Y \to Y'$. Since
$X'$ and $Y'$ form an admissible pair, Proposition~\ref{Padmissible} says that
if $F$ is not a homothety, then $F$ does not preserve some $\TP$ Toeplitz
kernel $K'$ on $X' \times Y'$. Transferring the kernel $K'$ back to a
totally positive kernel on $X \times Y$ via $( \varphi_X, \varphi_Y )$
shows the contrapositive of the non-trivial implication in
Theorem~\ref{Tgeneral}.

\subsection{Preservers of TP and TN on symmetric kernels}

We conclude by mentioning that there exist symmetrical analogues of all
of the results in this section, characterizing the preservers of symmetric
$\TN_p$ or $\TP_p$  kernels on $X \times X$ for an arbitrary totally
ordered set $X$. While the results have similar statements, we are unable
to employ P\'olya frequency functions or sequences as test kernels
in our proofs, and so are forced to look elsewhere.

Here we discuss only one result: the symmetric version of the final
theorem above.

\begin{theorem}[\cite{BGKP-TN}]
Given an infinite totally ordered set $X$ such that there exists a
symmetric $\TP_2$ kernel on $X \times X$, the only preservers of
symmetric $\TP$ kernels on $X \times X$ are the
positive homotheties.
\end{theorem}

The proof of this result employs all of the above tools and
techniques. We conclude this note with a sketch of it.

First, by the symmetric variant of Proposition~\ref{Pvdm}, every
preserver is a power function of the form
$F( x ) = c x^\alpha$ with $c > 0$ and $\alpha > 0$.
Moreover, as $F$ is continuous, it also preserves the class of
symmetric $\TN$ kernels on $X \times X$ that are limits of symmetric
$\TP$ kernels.

Using the order-preserving bijection $\varphi_X : X \to X'$ from
Proposition~\ref{Pembed}(3), we transfer a family of kernels
on $X' \times X'$, which are Hankel on $Z_n \times Z_n$ for each
$n$-step arithmetic progression $Z_n \subset X'$, to symmetric $\TN$
kernels on $X \times X$. As these Hankel kernels on $X' \times X'$ are
limits of $\TP$ Hankel kernels, the preceding paragraph and a variant
of the stronger Schoenberg-type Theorem~\ref{Thankel}
from~\cite{BGKP-hankel} imply that $F$ is absolutely monotonic. Hence
$\alpha \in \Z_{>0}$.

We now fix a real number $\beta > 0$ and consider the test function
\[
M_\beta : \R \to ( 0, \infty ); \ x \mapsto %
( \beta + 1 ) \exp( -\beta | x | ) - \beta \exp( -( \beta + 1 ) | x | ).
\]
It follows from Schoenberg's representation theorem,
Theorem~\ref{Tlaplace}, that $M_\beta$ is a $\PF$ function, but
that that $M_\beta^k$ is not a $\PF$ function for any integer $k \geq 2$.
Using convolution with a family of Gaussian densities, it is readily seen that
the Toeplitz kernel associated with $M_\beta$ is the limit of symmetric $\TP$ kernels,
and so $F \circ M_\beta = M_\beta^\alpha$ is $\TN$ on $X \times X$ by our
assumption on $F$.

Finally, we suppose $\alpha \geq 2$. Employing a discretization technique
and using the continuity of $M_\beta$, we produce an arithmetic progression
$Z' = ( z_1' < \cdots < z_{n_\alpha}' )$ in~$\R$ such that
$( M_\beta( z_i', z_j' )^\alpha )_{i, j = 1}^{n_\alpha}$ is not $\TN$.
We transfer this kernel on $Z' \times Z'$ first to $X' \times X'$, via a
change of scale and origin, then to $X \times X$ using $\varphi_X$.
This shows that $M_\beta^\alpha$ is not $\TN$ on $X \times X$, a contradiction, and
so the power~$\alpha$ cannot be $2$ or more, forcing $\alpha = 1$.

\subsection*{Acknowledgements}

D.G.~was partially supported by a University of Delaware Research
Foundation grant, by a Simons Foundation collaboration grant for
mathematicians, and by a University of Delaware strategic initiative
grant.

A.K.~was partially supported by Ramanujan Fellowship grant
SB/S2/RJN-121/2017 and SwarnaJayanti Fellowship grants SB/SJF/2019-20/14
and DST/SJF/MS/2019/3 from SERB and DST (Govt.~of India), and by grant
F.510/25/CAS-II/2018(SAP-I) from UGC (Govt.~of India).

M.P.~was partially supported by a Simons Foundation collaboration
grant.




\begin{thebibliography}{88}
\bibitem{ASW}
M.~Aissen, I.J.~Schoenberg, and A.M.~Whitney.
\newblock On the generating functions of totally positive sequences I.
\newblock \href{http://dx.doi.org/10.1007/BF02786970}%
{\em J.\ d'Analyse Math.}, 2:93--103, 1952.

\bibitem{Ando} 
T.~Ando.
\newblock Totally positive matrices.
\newblock \href{http://dx.doi.org/10.1016/0024-3795(87)90313-2}%
{\em Linear Algebra Appl.}, 90:165--219, 1987.

\bibitem{BCC}
A.~Bakan, T.~Craven, and G.~Csordas.
\newblock Interpolation and the Laguerre--P\'olya class.
\newblock {\em Southwest J.\ Pure Appl.\ Math.}, 1:38--53, 2001.

\bibitem{BGKP-TN}
A.~Belton, D.~Guillot, A.~Khare, and M.~Putinar.
\newblock Totally positive kernels, P\'olya frequency functions, and their transforms.
\newblock {\em J.\ d'Analyse Math.}, in press;
\href{http://arxiv.org/abs/2006.16213}{arXiv:math.FA/2006.16213}.

\bibitem{BGKP-fixeddim}
A.~Belton, D.~Guillot, A.~Khare, and M.~Putinar.
\newblock Matrix positivity preservers in fixed dimension. I.
\newblock \href{http://dx.doi.org/10.1016/j.aim.2016.04.016}%
{\em Adv.\ Math.}, 298:325--368, 2016.

\bibitem{BGKP-HW}
A.~Belton, D.~Guillot, A.~Khare, and M.~Putinar.
\newblock Hirschman--Widder densities.
\newblock \href{http://dx.doi.org/10.1016/j.acha.2022.04.002}%
{\em Appl.\ Comput.\ Harm.\ Anal.}, 60:396--425, 2022.

\bibitem{BGKP-hankel}
A.~Belton, D.~Guillot, A.~Khare, and M.~Putinar.
\newblock Moment-sequence transforms.
\newblock \href{http://dx.doi.org/10.4171/jems/1145}%
{\em J.\ Eur.\ Math.\ Soc.}, 24(9):3109--3160, 2022.

\bibitem{BFZ}
A.~Berenstein, S.~Fomin, and A.~Zelevinsky.
\newblock Parametrizations of canonical bases and totally positive matrices.
\newblock \href{http://dx.doi.org/10.1006/aima.1996.0057}%
{\em Adv.\ Math.}, 122:49--149, 1996.

\bibitem{BZ-2}
A.~Berenstein and A.~Zelevinsky.
\newblock Tensor product multiplicities, canonical bases and totally positive varieties.
\newblock \href{http://dx.doi.org/10.1007/s002220000102}%
{\em Invent.\ Math.}, 143:77--128, 2001.

\bibitem{Berezin}
F.A.~Berezin.
\newblock Quantization in complex symmetric spaces.
\newblock
\href{http://www.mathnet.ru/php/archive.phtml?wshow=paper&jrnid=im&paperid=1834&option_lang=eng}%
{\em Izv.\ Akad.\ Nauk SSSR Ser.\ Mat.}, 39(2):363--402, 1975.

\bibitem{Br1}
F.~Brenti.
\newblock {\em Unimodal, log-concave, and P\'olya frequency sequences in combinatorics}.
\newblock Mem.\ Amer.\ Math.\ Soc., vol.~413, American Mathematical Society, Providence, 1989.

\bibitem{Br2}
F.~Brenti.
\newblock Combinatorics and total positivity.
\newblock \href{http://dx.doi.org/10.1016/0097-3165(95)90000-4}%
{\em J.\ Combin.\ Theory Ser.~A}, 71(2):175--218, 1995.

\bibitem{Johnstone}
L.~Brown, I.M.~Johnstone, and K.B.~Macgibbon.
\newblock Variation diminishing convolution kernels associated with second-order differential operators.
\newblock \href{http://dx.doi.org/10.1080/01621459.1981.10477730}%
{\em J.\ Amer.\ Statist.\ Assoc.}, 76(376):824--832, 1981.

\bibitem{Projesh}
P.N.~Choudhury.
\newblock Characterizing total positivity: single vector tests via Linear Complementarity, sign non-reversal, and variation diminution.
\newblock {\em Preprint},
\href{http://arxiv.org/abs/2103.05624}{arXiv:math.RA/2103.05624}, 2021.

\bibitem{Curry2}
H.B.~Curry and I.J.~Schoenberg.
\newblock On P\'olya frequency functions IV: the fundamental spline functions and their limits.
\newblock \href{http://dx.doi.org/10.1007/BF02788653}%
{\em J.\ d'Analyse Math.}, 17:71--107, 1966.

\bibitem{Descartes}
R.~Descartes.
\newblock {\em Le G\'eom\'etrie}.
\newblock
\href{https://archive.org/details/geometryofrene00desc}%
{Appendix} to \textit{Discours de la m\'ethode}, 1637.

\bibitem{Edrei}
A.~Edrei.
\newblock On the generating functions of totally positive sequences II.
\newblock \href{http://dx.doi.org/10.1007/BF02786971}%
{\em J.\ d'Analyse Math.}, 2:104--109, 1952.

\bibitem{Efron}
B.~Efron.
\newblock Increasing properties of P\'{o}lya frequency functions.
\newblock \href{http://dx.doi.org/10.1214/aoms/1177700288}%
{\em Ann.\ Math.\ Stat.}, 36:272--279, 1965.

\bibitem{FJ}
S.~Fallat and C.R.~Johnson.
\newblock {\em Totally nonnegative matrices}.
\newblock Princeton Univ.\ Press, Princeton, 2011.

{\bibitem{Faraut}
J.~Faraut.
\newblock Rayleigh theorem, projection of orbital measures and spline functions.
\newblock \href{https://doi.org/10.1515/apam-2015-5012}%
{\em Adv.\ Pure Appl.\ Math.}, 6(4):261--283, 2015.}

\bibitem{FK-book}
J.~Faraut and \'A.~Kor\'anyi.
\newblock {\em Analysis on symmetric cones}.
\newblock
\href{https://global.oup.com/academic/product/analysis-on-symmetric-cones-9780198534778}%
{Oxford Mathematical Monographs}. The Clarendon Press, Oxford University Press, New York, second edition, 1994.
\newblock Oxford Science Publications.

\bibitem{Farrell}
R.H.~Farrell.
\newblock {\em Multivariate calculation: Use of the continuous groups}.
\newblock \href{https://doi.org/10.1007/978-1-4613-8528-8}%
{Springer Series in Statistics}, vol.\ 376, 1985.

\bibitem{Fe}
M.~Fekete and G.~P\'olya.
\newblock \"Uber ein Problem von Laguerre.
\newblock \href{http://dx.doi.org/10.1007/BF03015009}%
{\em Rend.\ Circ.\ Mat.\ Palermo}, 34:89--120, 1912.

\bibitem{FitzHorn}
C.H.~FitzGerald and R.A.~Horn.
\newblock On fractional {H}adamard powers of positive definite matrices.
\newblock \href{http://dx.doi.org/10.1016/0022-247X(77)90167-6}%
{\em J.\ Math.\ Anal.\ Appl.}, 61(3):633--642, 1977.

\bibitem{FZ-1}
S.~Fomin and A.~Zelevinsky.
\newblock Double Bruhat cells and total positivity.
\newblock \href{http://dx.doi.org/10.1090/S0894-0347-99-00295-7}%
{\em J.\ Amer.\ Math.\ Soc.}, 12:335--380, 1999.

\bibitem{FZ-2} 
S.~Fomin and A.~Zelevinsky.
\newblock Total positivity: tests and parametrizations.
\newblock \href{http://dx.doi.org/10.1007/BF03024444}%
{\em Math.\ Intelligencer}, 22(1):23--33, 2000.

\bibitem{GK}
F.R.~Gantmacher and M.G.~Krein.
\newblock Sur les matrices compl\`etement non n\'egatives et oscillatoires.
\newblock \href{http://www.numdam.org/item?id=CM_1937__4__445_0}%
{\em Compos.\ Math.}, 4:445--476, 1937.

\bibitem{GK1}
F.R.~Gantmacher and M.G.~Krein.
\newblock {\em Oscillation matrices and kernels and small vibrations of mechanical systems}.
\newblock Translated by A. Eremenko. Chelsea Publishing Co., New York, 2002.

\bibitem{TP}
M.~Gasca and C.A.~Micchelli.  
\newblock \href{http://dx.doi.org/10.1007/978-94-015-8674-0}%
{\em Total positivity and its applications}. 
Mathematics and its Applications, vol.~359,
\newblock Springer Science+Business Media, Dordrecht, 1996.

\bibitem{GN}
I.M.~Gel'fand and M.A.~Naimark.
\newblock Unitary representations of the classical groups.
\newblock \href{http://mi.mathnet.ru/eng/tm1100}%
{\em Trudy Mat.\ Inst.\ Steklova}, 36:3--288, 1950.

\bibitem{Gindikin}
S.G.~Gindikin.
\newblock Invariant generalized functions in homogeneous domains.
\newblock \href{http://dx.doi.org/10.1007/BF01078179}%
{\em Funct.\ Anal.\ Appl.}, 9:50--52, 1975.

\bibitem{Grochenig}
K.~Gr\"ochenig.
\newblock Schoenberg's theory of totally positive functions and the Riemann zeta function.
\newblock {\em Preprint},
\href{http://arxiv.org/abs/2007.12889}{arXiv:math.NT/2007.12889}, 2020.

\bibitem{GRS}
K.~Gr\"ochenig, J.L.~Romero, and J.~St\"ockler.
\newblock Sampling theorems for shift-invariant spaces, Gabor frames, and totally positive functions.
\newblock \href{http://dx.doi.org/10.1007/s00222-017-0760-2}%
{\em Invent.\ Math.}, 211:1119--1148, 2018.

\bibitem{GS}
K.~Gr\"ochenig and J.~St\"ockler.
\newblock Gabor frames and totally positive functions.
\newblock \href{http://dx.doi.org/10.1215/00127094-2141944}%
{\em Duke Math.\ J.}, 162(6):1003--1031, 2013.

\bibitem{Grommer}
J.~Grommer.
\newblock Ganze transzendente Funktionen mit lauter reellen Nullstellen.
\newblock \href{http://dx.doi.org/10.1515/crll.1914.144.114}%
{\em J. reine angew.\ Math.}, 144:114--166, 1914.

\bibitem{hc}
Harish-Chandra.
\newblock Differential operators on a semisimple Lie algebra.
\newblock \href{http://dx.doi.org/10.2307/2372387}%
{\em Amer.\ J.\ Math.}, 79(1):87--120, 1957.

\bibitem{HW49}
I.I.~Hirschman and D.V.~Widder.
\newblock The inversion of a general class of convolution transforms.
\newblock \href{http://dx.doi.org/10.1090/S0002-9947-1949-0032817-4}%
{\em Trans.\ Amer.\ Math.\ Soc.}, 66(1):135--201, 1949.

\bibitem{HW}
I.I.~Hirschman and D.V.~Widder
\newblock {\em The convolution transform.}
\newblock \href{https://press.princeton.edu/books/paperback/9780691626925/convolution-transform}%
{Princeton Legacy Library}, Princeton University Press, Princeton, 1955.

\bibitem{iz}
C.~Itzykson and J.-B.~Zuber.
\newblock The planar approximation. II.
\newblock \href{http://dx.doi.org/10.1063/1.524438}%
{\em J.\ Math.\ Phys.}, 21(3):411--421, 1980.

\bibitem{Jain2}
T.~Jain.
\newblock Hadamard powers of rank two, doubly nonnegative matrices.
\newblock \href{http://dx.doi.org/10.1007/s43036-020-00066-6}%
{\em Adv.\ in Oper.\ Theory}, 5:839--849, 2020 (Rajendra Bhatia volume).

\bibitem{James}
A.T.~James.
\newblock Distributions of matrix variates and latent roots derived from normal samples.
\newblock \href{http://dx.doi.org/10.1214/aoms/1177703550}%
{\em Ann.\ Math.\ Statist.}, 35(2):475--501, 1964.

\bibitem{KarlinTAMS}
S.~Karlin.
\newblock Total positivity, absorption probabilities and applications.
\newblock \href{http://dx.doi.org/10.1090/S0002-9947-1964-0168010-2}%
{\em Trans.\ Amer.\ Math.\ Soc.}, 111:33--107, 1964.

\bibitem{Karlin}
S.~Karlin.
\newblock {\em Total positivity. {V}olume 1.}
\newblock Stanford University Press, Stanford, 1968.

\bibitem{Katkova2007}
O.M.~Katkova.
\newblock Multiple positivity and the Riemann zeta-function.
\newblock \href{http://dx.doi.org/10.1007/BF03321628}%
{\em Comput.\ Methods Funct.\ Theory}, 7:13--31, 2007.

\bibitem{Khare2020}
A.~Khare.
\newblock Critical exponents for total positivity, individual kernel encoders, and the Jain--Karlin--Schoenberg kernel.
\newblock {\em Preprint},
\href{https://arxiv.org/abs/2008.05121}{arXiv:math.FA/2008.05121}, 2020.

\bibitem{Khare2022}
A.~Khare.
\newblock Smooth entrywise positivity preservers, a Horn--Loewner master theorem, and symmetric function identities.
\newblock \href{http://dx.doi.org/10.1090/tran/8563}%
{\em Trans.\ Amer.\ Math.\ Soc.}, 375(3):2217--2236, 2022.

\bibitem{KT}
A.~Khare and T.~Tao.
\newblock On the sign patterns of entrywise positivity preservers in
fixed dimension.
\newblock  \href{http://dx.doi.org/10.1353/ajm.2021.0049}%
{\em Amer.\ J.\ Math.}, 143(6):1863--1929, 2021.

\bibitem{KP}
J.S.~Kim and F.~Proschan.
\newblock
\href{http://dx.doi.org/10.1002/0471667196.ess2735.pub2}%
{Total positivity}. In: %
{\em Encyclopedia of Statistical Sciences} (S.~Kotz et al., Eds.),
vol.~14, pp.~8665--8672,
\newblock John Wiley \& Sons, New York, 2006.

\bibitem{KW-1}
Y.~Kodama and L.K.~Williams.
\newblock KP solitons, total positivity and cluster algebras.
\newblock \href{http://www.pnas.org/content/108/22/8984}%
{\em Proc.\ Natl.\ Acad.\ Sci.\ USA}, 108:8984--8989, 2011.

\bibitem{KW-2}
Y.~Kodama and L.K.~Williams.
\newblock KP solitons and total positivity for the Grassmannian.
\newblock \href{http://dx.doi.org/10.1007/s00222-014-0506-3}%
{\em Invent.\ Math.}, 198(3):637--699, 2014.

\bibitem{Laguerre1}
E.~Laguerre.
\newblock Sur les fonctions du genre z\'ero et du genre un.
\newblock {\em C.\ R.\ Acad.\ Sci.}, 95:828--831, 1882.

\bibitem{Laguerre2}
E.~Laguerre.
\newblock M\'emoire sur la th\'eorie des \'equations num\'eriques.
\newblock \href{http://sites.mathdoc.fr/JMPA/PDF/JMPA_1883_3_9_A5_0.pdf}%
{\em J.\ Math.\ Pures et Appl.}, 9:9--146, 1883.

\bibitem{LetacMassam}
G.~Letac and H.~Massam.
\newblock The Laplace transform $( \det s )^{-p} \exp {\rm tr}( s^{-1} w )$
and the existence of non-central Wishart distributions.
\newblock \href{http://dx.doi.org/10.1016/j.jmva.2017.10.005}%
{\em J.\ Multivar.\ Anal.}, 163:96--110, 2018.

\bibitem{Levin}
B.~Ja.~Levin.
\newblock {\em Distribution of zeros of entire functions.}
\newblock Translated from the Russian by R.P.\ Boas, J.M.\ Danskin, F.M.\ Goodspeed, J.\ Korevaar, A.L.\ Shields, and H.P.\ Thielman.
Revised edition. Translations of Mathematical Monographs~5.
American Mathematical Society, Providence, R.I., 1980.

\bibitem{Loewner55}
C.~Loewner.
\newblock On totally positive matrices.
\newblock \href{http://dx.doi.org/10.1007/BF01187945}%
{\em Math. Z.}, 63:338--340, 1955 (Issai Schur memorial volume).

\bibitem{L-1}
G.~Lusztig.
\newblock Total positivity in reductive groups, Lie theory and geometry.
\newblock {\em Progr.\ Math.}, vol.~123, Birkh\"auser Boston, Boston, MA, 1994, pp.531--568.

\bibitem{L-2}
G.~Lusztig.
\newblock Total positivity and canonical bases, Algebraic groups and Lie groups.
\newblock {\em Austral.\ Math.\ Soc.\ Lect.\ Ser.}, vol.~9, Cambridge 
University Press, Cambridge, 1997, pp.~281--295.

\bibitem{Mayer}
E.~Mayerhofer.
\newblock On Wishart and noncentral Wishart distributions on symmetric cones.
\newblock \href{http://dx.doi.org/10.1090/tran/7754}%
{\em Trans.\ Amer.\ Math.\ Soc.}, 371(10):7093--7109, 2019.

\bibitem{OV}
G.~Olshanski and A.M.~Vershik.
\newblock Ergodic unitarily invariant measures on the space of infinite Hermitian matrices.
\newblock In: \href{http://dx.doi.org/10.1090/trans2/175}%
{\em F.A.\ Berezin Memorial Volume}
(R.L.~Dobrushin, R.A.~Minlos, M.A.~Shubin, A.M.~Vershik, Eds.), 
Contemporary Mathematical Physics, Amer.\ Math.\ Soc.\ Transl.\ Ser.\ 2,
Vol.\ 175, Amer.\ Math.\ Soc., Providence, R.I., 1996, pp.137--175.

\bibitem{PR}
S.D.~Peddada and D.St.P.~Richards.
\newblock Proof of a conjecture of M.L.~Eaton on the characteristic function of the Wishart distribution.
\newblock \href{http://dx.doi.org/10.1214/aop/1176990455}%
{\em Ann.\ Probab.}, 19(2):868--874, 1991.
(See also \href{http://dx.doi.org/10.1214/aop/1176989822}%
{\em Ann.\ Probab.}, 20(2):1107, 1992.)

\bibitem{Pinkus}
A.~Pinkus.
\newblock {\em Totally positive matrices}.
\newblock Cambridge University Press, Cambridge, 2010.

\bibitem{Polya0}
G.~P\'olya.
\newblock \"Uber Ann\"aherung durch Polynome mit lauter reellen Wurzeln.
\newblock \href{https://doi.org/10.1007/BF03016033}%
{\em Rend.\ Circ.\ Mat.\ Palermo}, 36:279--295, 1913.

\bibitem{Postnikov}
A.~Postnikov.
\newblock Total positivity, Grassmannians, and networks.
\newblock {\em Preprint},
\href{http://arxiv.org/abs/math/0609764v1}{arXiv:math/0609764v1}, 2006.

\bibitem{Ri}
K.~Rietsch.
\newblock Totally positive Toeplitz matrices and quantum cohomology of partial flag varieties.
\newblock \href{http://dx.doi.org/10.1090/S0894-0347-02-00412-5}%
{\em J.\ Amer.\ Math.\ Soc.}, 16(2):363--392, 2003.

\bibitem{Rudin59}
W.~Rudin.
\newblock Positive definite sequences and absolutely monotonic functions.
\newblock \href{http://dx.doi.org/10.1215/S0012-7094-59-02659-6}%
{\em Duke Math.\ J.}, 26(4):617--622, 1959.

\bibitem{Schoenberg30}
I.J.~Schoenberg.
\newblock \"Uber variationsvermindernde lineare Transformationen.
\newblock \href{http://dx.doi.org/10.1007/BF01194637}%
{\em Math.\ Z.}, 32:321--328, 1930.

\bibitem{Schoenberg42}
I.J.~Schoenberg.
\newblock Positive definite functions on spheres.
\newblock \href{http://dx.doi.org/10.1215/S0012-7094-42-00908-6}%
{\em Duke Math.\ J.}, 9(1):96--108, 1942.

\bibitem{Schoenberg50}
I.J.~Schoenberg.
\newblock On P\'{o}lya frequency functions. II.~Variation-diminishing integral operators of the convolution type.
\newblock {\em Acta Sci.\ Math.\ (Szeged)}, 12:97--106, 1950.

\bibitem{Schoenberg51}
I.J.~Schoenberg.
\newblock On P\'{o}lya frequency functions. I.~The totally positive functions and their Laplace transforms.
\newblock \href{https://doi.org/10.1007/BF02790092}%
{\em J.\ d'Analyse Math.}, 1:331--374, 1951.

\bibitem{Schoenberg55}
I.J.~Schoenberg.
\newblock On the zeros of the generating functions of multiply positive sequences and functions.
\newblock \href{http://dx.doi.org/10.2307/1970073}%
{\em Ann.\ of Math.}, 62(3):447--471, 1955.

\bibitem{Schoenberg1973}
I.J.~Schoenberg.
\newblock {\em Cardinal spline interpolation.}
\newblock Conference Board of the Mathematical Sciences Regional Conference Series in Applied Mathematics,
no.~12. Society for Industrial and Applied Mathematics, Philadelphia, PA, 1973.

\bibitem{SchoenbergWhitney53}
I.J.~Schoenberg and A.M.~Whitney.
\newblock On P\'{o}lya frequency functions. III.
The positivity of translation determinants with an application to the interpolation problem by spline curves.
\newblock \href{https://doi.org/10.1090/S0002-9947-1953-0053177-X}%
{\em Trans.\ Amer.\ Math.\ Soc.}, 74:246--259, 1953.

\bibitem{Schur1911}
J.~Schur.
\newblock {B}emerkungen zur {T}heorie der beschr{\"a}nkten {B}ilinearformen mit unendlich vielen {V}er{\"a}nderlichen.
\newblock \href{http://dx.doi.org/10.1515/crll.1911.140.1}%
{\em J.\ reine angew.\ Math.}, 140:1--28, 1911.

\bibitem{Polya-Schur}
J.~Schur and G.~P\'olya.
\newblock \"Uber zwei Arten von Faktorenfolgen in der Theorie der algebraischen Gleichungen.
\newblock \href{https://doi.org/10.1515/crll.1914.144.89}%
{\em J.\ reine angew.\ Math.} 144:89--113, 1914.

\bibitem{Sra}
S.~Sra.
\newblock Positive definite functions of noncommuting contractions,
  {H}ua-{B}ellman matrices, and a new distance metric.
\newblock {\em Preprint},
\href{http://arxiv.org/abs/2112.00056}{arXiv:math.FA/2112.00056}, 2021.

\bibitem{Takemura}
A.~Takemura.
\newblock {\em Zonal polynomials}.
\newblock \href{https://www.jstor.org/stable/4355472}%
{IMS Lecture Notes Monograph Series}, vol.~4, 1984.

\bibitem{Thoma}
E.~Thoma.
\newblock Die unzerlegbaren, positiv-definiten Klassenfunktionen der abz\"ahlbar unendlichen, symmetrischen Gruppe.
\newblock \href{https://doi.org/10.1007/BF01114877}%
{\em Math.\ Z.}, 85:40--61, 1964.

\bibitem{Tschebotareff}
N.~Tschebotareff.
\newblock \"Uber die Realit\"at von Nullstellen ganzer transzendenter Funktionen.
\newblock \href{https://eudml.org/doc/159278}%
{\em Math.\ Ann.}, 99(1):660--686, 1928.

\bibitem{RossiVergne}
M.~Vergne and H.E.~Rossi.
\newblock Analytic continuation of the holomorphic discrete series of a semi-simple Lie group.
\newblock \href{http://dx.doi.org/10.1007/BF02392042}%
{\em Acta Math.}, 136:1--59, 1976.

\bibitem{VK}
A.M.~Vershik and S.V.~Kerov.
\newblock Characters and factor representations of the infinite unitary group.
\newblock \href{http://mi.mathnet.ru/eng/dan45717}%
{\em Sov.\ Math.\ Dokl.}, 26:570--574, 1982.

\bibitem{Voiculescu}
Dan Voiculescu.
\newblock Repr\'esentations factorielles de type $II_1$ de $U(\infty)$.
\newblock {\em J.\ Math.\ Pures et Appl.}, 55:1--20, 1976.

\bibitem{Wallach}
N.R.~Wallach.
\newblock The analytic continuation of the discrete series.
\newblock {\em Trans.\ Amer.\ Math.\ Soc.}, 251, 1979.
\href{http://dx.doi.org/10.1090/S0002-9947-1979-0531967-2}{Part I}, pp.~1--17;
\href{http://dx.doi.org/10.1090/S0002-9947-79-99965-3}{Part II}, pp.~19--37.

\bibitem{Whi}
A.M.~Whitney.
\newblock A reduction theorem for totally positive matrices.
\newblock \href{http://dx.doi.org/10.1007/BF02786969}%
{\em J.\ d'Analyse Math.}, 2(1):88--92, 1952.

\end{thebibliography}
\end{document}